\documentclass[12pt,leqno]{article}

 \usepackage{epigraph}
\usepackage{hyperref}
\usepackage{theoremref}
\usepackage{amsmath}
\usepackage{amsthm}
\usepackage{amssymb}
\usepackage{amsfonts}
\usepackage{color}
\usepackage[all]{xy} \xyoption{all}
 \usepackage{wrapfig}
\usepackage{enumerate}
\usepackage{enumitem}
\usepackage{mathabx}
\usepackage{url}

\usepackage{tikz, tikz-3dplot, pgfplots}
\usetikzlibrary{decorations.markings, patterns, decorations.pathmorphing}

\tikzset{->-/.style={decoration={
markings,
mark=at position #1 with {\arrow{>}}},postaction={decorate}}}

\def\CC{\mathbb{C}}

\def\GG{\mathbb{G}}

\def\RR{\mathbb{R}}
\def\ZZ{\mathbb{Z}}
\def\QQ{\mathbb{Q}}
\def\TT{{\mathbb{T}}}

\def\gen{\mathfrak{g}}

\def\len{\mathfrak{l}}

\def\sen{\mathfrak{s}}

\def\Hen{\mathfrak{H}}

\def\Ac{\mathcal{A}}

\def\Cc{\mathcal{C}}

\def\Dc{\mathcal{D}}
\def\Ec{\mathcal{E}}
\def\Fc{\mathcal{F}}
\def\Gc{\mathcal{G}}

\def\Lc{\mathcal{L}}
\def\Mc{\mathcal{M}}

\def\Hc{\mathcal{H}}
\def\Oc{\mathcal{O}}

\def\Qc{\mathcal{Q}}

\def\Sc{\mathcal{S}}

\def\Xc{\mathcal{X}}

\def\frakg{\mathfrak{g}}

 \def\bPhi{{\mathbf{\Phi}}}
\def\bPsi{{\mathbf{\Psi}}}
 \def\bv{{{\bf v}}}

 \def\be{\begin{equation}}
 \def\bef{\begin{figure}}
 \def\bem{\begin{matrix}}
 \def\beps{{\mathbf{\varepsilon}}}
 
 \def\bpm{\begin{pmatrix}}
 
 \def\btp{\begin{tikzpicture}}

 \def\centerarc[#1](#2)(#3:#4:#5)
{ \draw[#1] ($(#2)+({#5*cos(#3)},{#5*sin(#3)})$) arc (#3:#4:#5); }
 
 \def\ch{{\on{char}}}

 \def\Crit{{\on{Crit}}}
 \def\CS{{\on{CS}}}

 \def\del{{\partial}}
 \def\dirr{{\on{dir}}}

  \def\ee{\end{equation}}
 
 \def\enf{\end{figure}}
 \def\enm{\end{matrix}}
 \def\epm{\end{pmatrix}}
 \def\eps{{\varepsilon}}
 \def\etp{\end{tikzpicture}}

 \def\fr{{\on{fr}}}
 \def\FT{{\on{FT}}}

 \def\gen{{\on{gen}}}

 \def\hol{{\on{h}}}
 \def\hr{{\on{hr}}}
  \def\Hom{\operatorname{Hom}\nolimits}

   \def\hra{\hookrightarrow}

 \def\Id{\operatorname{Id}\nolimits}

  \def\Im{{\on{Im}}}

 \def\k{\mathbf k}
  \def\Ker{\operatorname{Ker}\nolimits}

 \def\Li{{\on{Li}}}

   \def\lra{\longrightarrow}
\def\LS{\on{{LS}}}

  \def\Mod{{\text{-}\on{Mod}}}

\def\Ob{\operatorname{Ob}\nolimits}
\def\on{\operatorname}
\def\ol{\overline}
\def\oo{{\infty}}

 \def\Perv{{\on{Perv}}}
 \def\perv{{\on{perv}}}
 \def\phi{{\varphi}}

 \def\rk{{\on{rk}}}
\def\Re{{\on{Re}}}
 \def\rect{{\on{rect}}} 
 \def\Rep{{\on{Rep}}}
 \def\RHom{{\on{RHom}}}

  \def\Sol{{\on{Sol}}}
  \def\Spec{{\on{Spec}}}
  \def\St{{\on{St}}}

\def\tr{{\on{tr}}}

\def\ul{\underline}

\def\wh{ \widehat}
\def\wt{\widetilde}

\def\1{{\mathbf{1}}}
\def\2{{\mathbf{2}}}
\def\(({(\hskip -1mm (}
\def\)){)\hskip -1mm )}
\def\-{{\setminus}}
\def \= { {\, \simeq \,}} 
\def\be{\begin{equation}}
\def\ee{\end{equation}}
\def\ed{\end{document}}

\def\la{{\langle}}
\def\ra{{\rangle}}

\numberwithin{equation}{subsection}

\newtheorem{thm}[equation]{Theorem}
\newtheorem{cor}[equation]{Corollary}
\newtheorem{lem}[equation]{Lemma}
\newtheorem{prop}[equation]{Proposition}
\newtheorem{refo}[equation]{Reformulation}

\theoremstyle{definition}
\newtheorem{definition}[equation]{Definition}
\newtheorem{defi}[equation]{Definition}
 
\theoremstyle{remark}
\newtheorem{remark}[equation]{Remark}
\newtheorem{rem}[equation]{Remark}
\newtheorem{rems}[equation]{Remarks}

\newtheorem{ex}[equation]{Example}
\newtheorem{exas}[equation]{Examples}

\numberwithin{itemcounter}{subsection}
\numberwithin{equation}{subsection}
\usepackage[toc,page]{appendix}

\usepackage{etoolbox}
\appto\appendix{\addtocontents{toc}{\protect\setcounter{tocdepth}{1}}}

\begin{document}

\title{ Resurgence  and perverse sheaves}

\author{Mikhail Kapranov,  Yan Soibelman}

	\maketitle

 {\em To our friend Maxim Kontsevich on occasion of his sixtieth birthday } 
 
 \vskip 1cm
 
 	\begin{abstract}
	
 We propose a point of view on resurgence theory based on the  study of
 perverse sheaves on the complex line  carrying an algebraic structure
 with respect to additive convolution. 
 In particular, we lift the  concept of alien derivatives
 introduced originally by J. \'Ecalle, 
  to the framework of perverse sheaves and study its behavior under sheaf-theoretic
  convolution. The full fledged
 resurgence theory needs a (yet undeveloped) generalization of the concept of perverse sheaves  allowing infinite, possibly dense,  sets of singularities.
  We discuss possible approaches
 to defining such objects and  some potential
 examples of them  coming from Cohomological Hall Algebras,
 wall-crossing structures and Chern-Simons theory. 

 \end{abstract}


 \addtocounter{section}{-1}
 
 \tableofcontents

 
 \epigraph{Le sommeil est plein de miracles!
 \\
Par un caprice singulier,
\\
J’avais banni de ces spectacles
\\
Le v\'eg\'etal irr\'egulier.}
{Ch. Baudelaire, {\em R\^eve parisien}}

 \section{Introduction} 
 
 The theory of resurgent functions pioneered by J. \'Ecalle \cite{ecalle}  studies analytic functions
 ``given by divergent series'' in terms of singularities of their Borel (formal Fourier-Laplace)
 transform (see \S \ref{subset:bor-bor} below for discussion). As with  any kind of Fourier transform, 
 this procedure takes multiplication into additive convolution. Resurgent analysis then proceeds
 by using equations that involve convolution as well as the monodromy data
 (alien derivatives) for the Borel-transformed functions.
 
 \vskip .2cm
 
 The goal of this paper is to propose a more conceptual point of view on resurgence theory
 by using the notion of perverse sheaf. It is known that operations such as convolution
 or Fourier transform can be defined at the sheaf-theoretical level and match, to some extent,
 the corresponding analytic operations for functions when such functions are realized as
 sections of the sheaves in question.  So, it is a natural idea to import the
 perverse sheaves language  into resurgent analysis, with the hope to achieve
 greater conceptual clarity.
 
 \vskip .2cm
 
 In a few words, our proposal is to consider perverse sheaves on $\CC$ carrying the
 structure of an algebra with respect to convolution (see \S \ref{subsec:gen-prog}
 for more details). Using this structure one can write resurgent equations whose
 unknowns will be sections of such ``resurgent perverse sheaves'', with classical examples
 appearing when sections are realized as actual analytic functions. But one can also,
 in principle, construct resurgent perverse sheaves in a more abstract fashion,
 not unlike the way one constructs field extensions not necessarily embedded in $\CC$
 by adjoining roots of algebraic equations. 
 
 \vskip .2cm
 
 Already in 1985 B. Malgrange  \cite{malgrange-ecalle} gave an interpretation of the resurgent
 formalism in terms of M. Sato's theory of microfunctions.  His main observation was that
 the concept of  a ``singularity'' (or ``singular part''), ubiquitous in this formalism, is but
 a synonym for a microfunction.  Now, the space of vanishing cycles of a perverse sheaf on
 $\CC$ is the same as the space of microfunction solutions of the corresponding
 holonomic regular  D-module \cite{GGM}. 
 Therefore, working with perverse sheaves and their vanishing cycles
 is a natural conceptual framework for resurgence theory. 
 
 \vskip .2cm
 
 However, for true applications to resurgence the theory of perverse sheaves must be extended
 to match the kind of multivaluedness that resurgent  functions typically possess. These functions
 typically have infinite  or even dense sets of singular points.  
 That is,  on any ``branch'' the singular points
 are of course discrete, but  going around each one
 leads to a new branch with new singularities etc.  Such behavior is referred to as
 ```analytic continuation without end''. 
 
  \vskip .2cm 
 
 In this paper we  do not attempt to generalize
 the theory of perverse sheaves in this direction. Instead, we  develop a resurgence-like formalism
 involving the standard concept of perverse sheaves on $\CC$ (i.e.  with finitely many
 singularities). Already this allows us to highlight many of the familiar  features  in the
 sheaf-theoretic context, for example the interpretation of  Stokes data via
 Picard-Lefschetz type formulas in the Borel plane. 
 Further, some examples of ``resurgent perverse sheaves'' may  be  already
 given   in this restricted context,  such as the version of Cohomological Hall algebra
 in \S \ref {subsec:COHA}. 
 
 \vskip .2cm
 
 A special role in our considerations is played by the category $\ol\Perv(\CC)$
 obtained by localizing $\Perv(\CC)$,  the abelian category of all perverse sheaves
 on $\CC$, by the Serre subcategory of (shifted) constant sheaves \cite {gelfand-MV}. 
 Objects of $\ol\Perv(\CC)$ have well-defined ``tunnelling data'' consisting of
 the spaces of vanishing cycles $\Phi_a, a\in\CC$ and the transport maps $m_{ab}(\gamma): \Phi_a\to\Phi_b$ for various paths $\gamma$ joining $a$ and $b$.

  \vskip .2cm
 It is known
 \cite{KKP, fressan-jossen} that $\ol\Perv(\CC)$ can be embedded back into $\Perv(\CC)$
 as the subcategory $\Perv^0(\CC)$ formed by perverse sheaves $\Fc$ with
 $H^\bullet(\CC,\Fc)=0$. 
 The operation of additive convolution $\Fc *\Gc$  is most easily defined using this
 realization \cite{fressan-jossen}. Our ``toy resurgent formalism'' can be
 seen as further  study of the Tannakian Galois group of the tensor category
 $(\ol\Perv(\CC), *)$, as defined and already studied in  \cite{fressan-jossen}. 
 So, our larger point is that the full fledged resurgent formalism is just a similar
 study but for a more general concept of perverse sheaves, still to be
 defined rigorously (see \S \ref{subsec:gen-prog} below). 
 We plan to discuss this  further in a future work.

 \vskip .2cm
 
 The paper consists of three chapters. In Chapter \ref {sec:perv-four} we
 recall  the elementary theory of perverse sheaves on $\CC$ (with finitely many singularities)
 with emphasis on features that we need.  The motivational \S \ref{subset:bor-bor} 
 explains the general framework of Borel summation and the resulting ``doctrine of
 two planes'': the original one carrying  functions given by
 divergent series and the Borel one where
 things become more topological.  In \S \ref{subsec:perRiem} we recall the
 basic definitions and emphasize the Picard-Lefschetz formula (Proposition
 \ref {prop:PL-form}) in the context
 of perverse sheaves. In
 \S \ref {subsec:locper} we present the Gelfand-MacPherson-Vilonen classification of
 perverse sheaves and of objects of the localized category $\ol\Perv(\CC, A)$, where
 constant sheaves are factored out but the vanishing cycles and transport maps remain. 
 The Fourier transform for perverse sheaves is explained in
 \S  \ref{subsec:FT}. There, we mostly follow \cite{KSS}. In \S \ref {subsec:lefper}
 we discuss a particular class of examples of perverse sheaves on $\CC$  associated to a regular function
 $S: X\to\CC$ on a complex algebraic variety. We call them
 {\em Lefschetz perverse sheaves}  $\Lc_S$.   The Fourier transform of $\Lc_S$
 can be seen as a categorification of the exponential integral associated to $S$. 
 
 \vskip .2cm 
 
 In Chapter \ref {sec:alien} we build up features of perverse sheaves on $\CC$
 which are most remindful of resurgence formalism. In this, we extend the
 analysis of  rectilinear transports given in \cite{KSS}, in  the case when the set $A$  of
 singularities is in linearly general position, to the arbitrary case, when an interval
 $[a,b], a, b\in A$ can contain intermediate points. 
  In \S \ref {subsec:conv} we discuss additive convolution of perverse sheaves on $\CC$. 
  Motivated by the classical resurgent formalism,  we  consider,  in \S \ref {subsec:trans-avoid},
  various ways of modifying the rectilinear transport so as to avoid  the intermediate
  points. The formulas for alien derivatives appear naturally in this context
  as some linear combinations of such  modified transports. In \S \ref {subsec:alien-der}
  we study alien derivatives more systematically; we also explain their
  relation with the Stokes automorphisms for the Fourier transform. 
  
  \vskip .2cm
  
  In the final Chapter \ref {sec:resurg} we discuss how our approach can be applied
  to actual resurgence problems.  This chapter is more speculative. 
  We start by sketching in
 \S \ref {subsec:gen-prog} the general program of studying perverse sheaves which
 are algebras with respect to the convolution,  highlighting the difficulties
 that are present in the general case. Then, we discuss several classes of
 potential examples: the example with  the Cohomological Hall algebra (a.k.a COHA) of a quiver in \S \ref {subsec:COHA},
 that of ``cluster perverse sheaves'' associated to  wall-crossing structures in \S 
 \ref {subsec:cluster}, and that of  Lefschetz perverse sheaves associated to the 
 complex Chern-Simons functional  in \S \ref {subsec: CS}. 
 
 \vskip .5cm
 
 {\it Acknowledgments.} The research of M.K. was supported by World Premier International Research Center Initiative
  (WPI Initiative),
 MEXT, Japan
 and  by the JSPS  KAKENHI grant 20H01794. 
 The work of Y.S. was partially supported by Simons Foundation grant MP-TSM-00001658, NSF grant 2401518, IHES Sabbatical Professorship and ERC grant "Renew Quantum".
 
 We would like to express our gratitude to A. Soibelman for correcting the English in the Introduction.

  
 \section {Perverse sheaves and their Fourier transform} \label{sec:perv-four}

 \subsection{Motivation: Borel summation and the Borel plane} \label{subset:bor-bor}
 
  \numberwithin{equation}{subsection}
 
 The famous Borel summation process for divergent series can be seen as an
 application of the Fourier transform in the complex domain.   It connects
 two copies of the complex plane $\CC$ which are loosely related 
 ``by the Fourier transform'': 
 
 \begin{figure}[h]
 \centering
 \begin{tikzpicture}[scale=0.6]
 \draw (-4,-4) -- (-4,4) -- (4,4) -- (4, -4) -- (-4, -4); 
  \draw (3,4) -- (3,3) -- (4,3); 
 
 \node at (3.5,3.5) {\large$w$}; 
 
  \node at (1,2 ){\small$\bullet$}; 
   \node at (3,0.3){\small$\bullet$}; 
      \node at (-2,1.5){\small$\bullet$}; 
      \node at (-1,2.5){\small$\bullet$};    
    \node at (-2,-2){\small$\bullet$};    
      \node at (3,-1.5){$\small\bullet$};     
\node at   (0, -2.5){$\cdots$};  
\node at (0,0){(Singularities)} ;  
 \end{tikzpicture}
 \begin{tikzpicture}[scale=0.6]
 \draw[color={white}]    (-4,-4) -- (-4,4) -- (4,4) -- (4, -4) -- (-4, -4);
 \draw   [ decoration={markings,mark=at position 0.99 with
{\arrow[scale=1.5,>=stealth]{>}}, 
mark=at position -1with
{\arrow[scale=1.5, >=stealth]{<}}
},postaction={decorate},
line width=.4mm] (-3.5,0) -- (3.5,0);
\node at (0, 0.5){Fourier}; 
\node at (0, -0.5){transform}; 
 \end{tikzpicture}
 \begin{tikzpicture}[scale=0.6]
 \draw (-4,-4) -- (-4,4) -- (4,4) -- (4, -4) -- (-4, -4);
 \draw (3,4) -- (3,3) -- (4,3); 
  \node at (3.5,3.5) {\large$z$}; 
  \node at (0,0){$\bullet$}; 
 \node at (-0.5, -0.5){$0$}; 
 \filldraw[opacity = 0.2] (0,0) -- (-2,4) -- (2, 4) -- (0,0); 
 \draw (-2,4) -- (0,0) -- (2,4); 
 \node at (0, -1.5){(Sectors)}; 
 \end{tikzpicture}
 \end{figure}
 
 \begin{itemize}
\item   The original (``irregular'') plane
 $\CC_z$ with coordinate (``large parameter'') $z$ in which we study possibly divergent
 formal power series near $\oo$:
 \begin{equation}\label{eq:f-hat}
 \wh f(z) = \sum_{n=0}^\oo {a_n\over z^{n+1}}
 \end{equation}
 Such series typically satisfy linear differential equations irregular at $\oo$,
 are divergent everywhere but serve as asymptotic expansions of interesting
 analytic solutions in various sectors. Quantum mechanical asymptotic series in powers of the
 Planck constant $\hbar$ are realized here by putting $z=1/\hbar$. 
 
 \item The  dual (Borel or ``regular'') plane $\CC_w$ with coordinate $w$ where we study
 (or obtain) solutions of differential equations with regular singularities, given by convergent
 series in $w$, which represent multivalued functions and  so can be analyzed
 topologically. The modern way of doing so is by using the language of perverse
 sheaves. 
 \end{itemize} 
 
 \noindent On the formal level, Borel summation of the series \eqref{eq:f-hat}
 (when it is possible) consists, first, of taking the termwise Fourier transform
 of the series using the identity  (particular case
  of the general formula  for the Fourier transform 
 of $z^\alpha$)
 \[
 \FT \biggl( {1\over z^{n+1}}\biggr) = {w^n\over{n!}}. 
 \]
 This gives the Borel-transformed series
 \[
 \wh f^B(w) = \sum_{n=0}^\oo {a_n\over{n!}} w^n
 \]
 which has more chances to converge. In good cases it has nonzero radius of
 convergence and extends to an analytic
 function $f^B(w)$  ``on  the entire $\CC_w$''  but possibly multivalued, with singularities etc. 
 Then the prescription for the sum $f(z)$ of the original series $\wh f(z)$
 is obtained by taking the inverse\footnote{the kind that sends $w^n/n!$
 back into $1/z^{n+1}$} Fourier (-Laplace) transform of $\wh f^B(w)$:
 \[
 f(z) := {\int_0^\oo}  f^B(w) e^{-zw} dw. 
 \]
Because of singularities of $f^B(w)$ there can be several inequivalent allowable
choices of the integration contour leading to ambiguity of the Borel sum,
which is not surprising if the series is divergent. 

\vskip .2cm

On the conceptual level, the Borel summation approach can be said to consist in 
representing irregular data\footnote{ Here we understand the word ``irregular\rq\rq{} in the wider sense, including but not restricted to differential equations and functions
satisfying them.} as Fourier transforms of regular ones. The success of this approach
comes from the fact that  such representation is possible in  many cases of
practical interest. Differential equations which are Fourier transforms
of regular ones form a rather special class.

 \vfill\eject

 \subsection{Perverse sheaves on  Riemann surfaces }\label{subsec:perRiem}
 
 \paragraph{Generalities.}
 We fix a base field $\k$.  Let $X$ be a complex analytic manifold.
 We denote by $\LS(X)$ the category of local systems of
 (not necessarily finite-dimensional) $\k$-vector spaces on $X$.
 
 Let $\Sc$ a
 locally finite complex Whitney stratification of  $X$. Thus each stratum $S\in\Sc$
 is a complex submanifold; we denote $i_S: S\to X$ the embedding.
 The closure $\ol S$ is, in general, a singular complex space. 
 We denote by $D^b(X,\Sc)$ the triangulated category of
  bounded  $\Sc$-constructible
  complexes of sheaves $\Fc$ on $X$. By definition, $\Sc$-constructibility
 of $\Fc$ means that
 each cohomology sheaf $\ul H^i(\Fc)$ is $\Sc$-constructible.
  That  is, each $i_S^* \ul H^q(\Fc)= \ul H^q(i_S^*\Fc)$   is an object of $\LS(S)$.

 We denote $\Perv(X, \Sc)\subset D^b(X,\Sc)$ the abelian category of  
 perverse sheaves of $\k$-vector
 spaces on $X$ smooth with respect to $\Sc$.  Explicitly, $\Perv(X,\Sc)$
 consists of complexes $\Fc\in D^b(X,\Sc)$   with the following properties:
 \begin{itemize}
 
  \item[(1)] For each stratum $S\in\Sc_A$ we have $\ul H^q(i^*_S\Fc) = 0$ for
 $q> -\dim_\CC S$.
 
 \item[(2)]  For each stratum $S\in\Sc_A$ we have $\ul H^q(i^!_S\Fc)=0$ for
 $q<-\dim_\CC S$. 
 \end{itemize}
 
 \noindent For example, if $\Lc\in \LS(X)$ is a local system, then $\Lc[\dim X]$,
 i.e., $\Lc$ put in degree $(-\dim X)$, is perverse and   lies in $\Perv(X,\emptyset)$. 
 
 \begin{rem}
 In \cite{KSS} we used a different normalization of the perversity conditions
 for which a local system in degree $0$ is considered to be perverse. 
 In  references to \cite{KSS} later in this paper, this difference, being
 easy to account for, is not further highlighted. 
 
 \end{rem}
 

 \paragraph{Nearby and vanishing cycles.} \label{par:nearvan}
  From now on we assume $\dim_\CC X = 1$, so $X$ is a complex curve (Riemann surface,
 possibly non-compact). Then a stratification $\Sc$ of $X$ is given by
 a discrete subset $A\subset X$ so the strata are elements of $A$ and the
 complement $X\- A$. In this case we use the notation  $\Perv(X,A)$ for $\Perv(X,\Sc)$. 
 The definition implies that for  any $\Fc\in\Perv(X,A)$ the restriction
 $\Fc|_{C\- A}$ is quasi-isomorphic to a local system
 (not necessarily of finite rank) in degree $-1$ and that $\ul H^q(\Fc) =0$ for $q\neq -1, 0$.

 \begin{ex}\label{ex:Phi-Psi}
 Let $X=D = \{|z|<1\}$ be the unit disk in $\CC$ and $A=\{0\}$. 
 In this case it is classical  \cite{GGM}  that $\Perv(D,0)$ is  equivalent to the category  of diagrams
\[
\xymatrix{
\Phi \ar@<.4ex>[r]^u&\Psi \ar@<.4ex>[l]^v
}
\]
of  $\k$-vector spaces
are linear maps such that
$T_\Psi: = \Id_\Psi - uv$   is an isomorphism
(or, what is equivalent, such that  $T_\Phi=\Id_\Phi-vu$ is an isomorphism). 
\end{ex}

More precisely, see  \cite{GGM} \cite[Prop. 1.1.6]{KSS}, the equivalence above depends
on a choice of a {\em radial cut}  $K\subset D$, a simple curve starting at 0 and ending on
the boundary $\del D$.  
\begin{figure}[h]
\centering
\begin{tikzpicture}[scale =0.4]

\node at (0,0) {\small$\bullet$};
\draw[line width=0.5] (0,0) circle (3cm); 
\draw[line width = 0.4mm]  plot [smooth] coordinates {(0,0)  (1, 0.1)  (2, -0.1) (3,0)}; 
\node at (-0.5, -0.5){$0$}; 
\node at (3, 3){$D$}; 
\node at (2, 0.6){$K$}; 

\end{tikzpicture}
\end{figure}
Given such $K$,  the sheaves of $K$-supported hypercohomology $\ul H^q_K(\Fc)$ are $0$ for $q\neq 0$,  the   sheaf
$\ul H^0_K(\Fc)$ on $K$ is constructible with respect to the stratification of $K$ into
$\{0\}$ and $K-\{0\}$, the spaces $\Phi$ and $\Psi$ associated to $\Fc$ are found as its stalks:
\[
\Phi = (\ul H^0_K(\Fc))_0, \quad \Psi =   (\ul H^0_K(\Fc))_\eps, \,\,\forall \eps\in K\-\{0\},
\]
and the map $u: \Phi\to\Psi$ is the generalization map of $\ul H^0_K(\Fc)$. See
 \cite{GGM} and the discussion after \cite[Prop. 1.1.6]{KSS} for the definition of
 $v: \Psi\to\Phi$.  The spaces $\Phi$ and $\Psi$ are called the spaces
 of {\em vanishing cyclies} and {\em nearby cycles} of $\Fc$ at $0$
 (in the direction of $K$). Note also that $\Psi$ is identified with the stalk of 
 the local system $(\Fc[-1])|_{D\-\{0\}}$ at any $\eps\in K-\{0\}$
 (hence the name ``nearby cycles''). 
 
 \vskip .2cm
 
 More generally, for a Riemann surface $X$ and a point $a\in X$, we denote
 by $S^1_a = S^1_a(X)$ the circle of directions at $a$. If $K$ is a smooth
 simple curve ending at $a$, we denote by $\dirr_a(K)\in S^1_a$ the direction
 of $K$ at $a$. If $A\subset X$ is discrete and $\Fc\in\Perv(X,A)$, then in a small disk
 near any $a\in A$ we have the situation of Example \ref{ex:Phi-Psi}.
 In particular, the vector spaces of vanishing and nearby cycles of $\Fc$ at
 $a$, being dependent on the direction of a cut, are naturally local systems on $S^1_w$
 which we denote $\bPhi_w(\Fc)$ and $\bPsi_w(\Fc)$. We thus have
 functors
 \be\label{eq:PhiPsi-LS}
 \bPhi_a, \bPsi_a: \Perv(X,A) \lra \LS(S^1_a).
 \ee

 \paragraph{Transport maps.} \label{par:transport}
We now recall the construction of {\em curvilinear transport maps} from
 \cite[\S 1.1C]{KSS}
 Let $\dim_\CC(X)=1$ and $A\subset X$ be discrete.
 Let $\alpha$ be a simple, piecewise smooth  arc in $X$ joining 
 two distinct   points $a,b\in A$
and not passing through any other elements of $A$  see  Fig. \ref {fig:transport}.
 Let us  equip $\alpha$ with the orientation going from $a$ to $b$.

 Considering $\alpha$ as a closed subset in $X$,
we have the sheaf $\ul H^0_\alpha(\Fc)$ on $\alpha$ which is constant on the open arc 
$\alpha-\{a,b\}$.

 \vskip 1cm
 
 \bef[h]
\centering
\btp [scale=.8, baseline=(current  bounding  box.center)]

\node (i) at (0,0){};
\fill (i) circle (0.1);

\node (j) at (12,0){};
\fill (j) circle (0.1);

\node (k) at (5,-2){};
\fill (k) circle (0.1);

\node (b) at (4,0.3){$$};
\fill (b) circle (0.08);

\draw [dotted, line width =0.3mm]  (0,0) circle (1cm);
\draw [dotted, line width =0.3mm]  (12,0) circle (1cm);

\draw[->, line width = 0.2mm] plot [smooth,tension=1.5] coordinates{
(0,0) (4,0.3)
(8,-0.3)
(12,0)
};
\node at (-0.5,0){$a$};
\node at (12.5, 0){$b$};
\node at (5.5, -2){$c$};
\node at (9,-0.7){\large$\alpha$};
\node at (4,-0.3){$p$};
\node at (4,0.8){$\Psi_\alpha$};
\node at (1.4, 0.7){$\Phi_{a,\alpha}$};
\node at (10.5, 0.3){$\Phi_{b,\alpha}$};

\etp

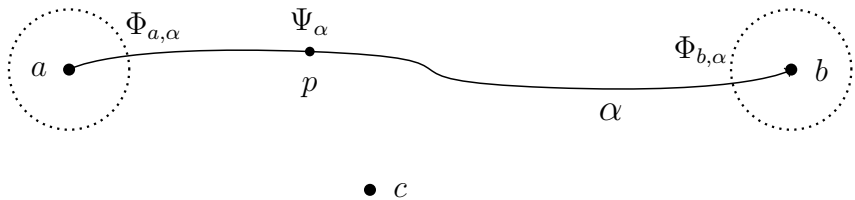
\caption{ The transport map.  }
\label{fig:transport}
\enf
Its stalks at $a$ and $b$ are the vanishing cycle spaces $\Phi_{a,\alpha}$,
$\Phi_{b,\alpha}$ for $\Fc$ at $a$ or $b$ in the direction of $\alpha$,
while its restriction to $\alpha\-\{a,b\}$ is a local system on the open interval;
in particular, all the stalks of this local system are canonically identified;
let us denote their common value by $\Psi_\alpha$. 
So we have the maps
\be\label{eq:u-a,alpha}
\xymatrix{
\Phi_{a,\alpha}  \ar@<.4ex>[r]^{u_{a,\alpha}}&\Psi_\alpha  \ar@<.4ex>[l]^{v_{a,\alpha}}
\ar@<-.4ex>[r]_{v_{b,\alpha}}
&
\Phi_{b, \alpha}
\ar@<-.4ex>[l]_{u_{b,\alpha}},
}
\ee
obtained from the description of $\Fc$ on small disks near $a$ and $b$ and using 
$\gamma$ as  the choice for a cut $K$.
We define the {\em transport map}
along $\alpha$ as
\be\label{eq:transport}
m_{ab}(\alpha) =m_{ab}^\Fc(\alpha) \,:= \, v_{b,\alpha}\circ u_{a, \alpha}: \Phi_{a,\alpha} \lra \Phi_{b,\alpha}.
\ee

\paragraph{Picard-Lefschetz identities for transports.} 
The construction of the maps $m_{ab}(\alpha)$  being purely
topological,  it is unchanged 
under isotopic deformations of $\alpha$ which do not pass through other 
elements of
$A$. More precisely   \cite[\S 1.1C]{KSS}, let 
$(\alpha_t)_{t\in[0,1]}$ be an 
  {\em admissible isotopy} of paths from $a$ to $b$, i.e., 
   a  continuous $1$-parameter family of simple arcs
$(\alpha_t)_{t\in[0,1]}$,
each $\alpha_t$ joining $a$ with $b$ and not passing through any other $c\in A$.
Then   we have a commutative
diagram
\be\label{eq:def-inv}
\xymatrix{
\Phi_{a,\alpha_0}
\ar[d]_{t_a}
\ar[r]^{m_{ab}(\alpha_0)} & \Phi_{b, \alpha_0}\ar[d]^{t_b}
\\
\Phi_{a,\alpha_1} \ar[r]^{m_{ab}(\alpha_1)} & \Phi_{b, \alpha_1},
}
\ee
where $t_a$ is the monodromy of the local system $\bPhi_a$ on $S^1_{w_i}$ from $\dirr_a(\alpha_0)$ to
$\dirr_a(\alpha_1)$, and similarly for $t_b$.

We now recall what happens when a path crosses a single point of $A$.
  That is, we consider a situation as in Fig. \ref{fig:PicLef},
where a path $\gamma'$ from $a$ to $c$ approaches the composite path formed by $\beta$ from $a$ to $b$
and $\alpha$ from $b$ to $c$. After crossing $b$, the path $\gamma'$ is changed to 
$\gamma$.

\bef[h]
\centering
\btp[scale=.8, baseline=(current  bounding  box.center)]

\node (i) at (3.5, -3){};
\fill(i) circle (0.15);

\node (k) at (-3, -3.5){};
\fill (k) circle (0.15);

\node (j) at (.7, -1){};
\fill (j) circle (0.15);

 \centerarc[line width=0.5](0.7, -1)(325:580:0.7)

\draw [ decoration={markings,mark=at position 0.7 with
{\arrow[scale=1.5,>=stealth]{>}}},postaction={decorate},
line width = .3mm]   (-3,-3.5) .. controls (1,2) ..  (3.5, -3) ;

\draw [decoration={markings,mark=at position 0.7 with
{\arrow[scale=1.5,>=stealth]{>}}},postaction={decorate},
line width = .3mm]   (-3,-3.5) .. controls (1,-4) ..   (3.5,-3);

\draw [decoration={markings,mark=at position 0.7 with
{\arrow[scale=1.5,>=stealth]{>}}},postaction={decorate},
line width = .3mm]  (0.7,-1)  .. controls (2,-2) ..  (3.5,-3) ;

\draw [ decoration={markings,mark=at position 0.7 with
{\arrow[scale=1.5,>=stealth]{>}}},postaction={decorate}, line width = .3mm]  (-3,-3.5) .. controls (-2,-3) .. (0.7,-1);

\node at (4.5, -3.5){$c$};

\node at ( -4,-3.5){$a$};
\node at ( 0.7,-0.6){$b$};
\node at (0,-4.2){\large$\gamma'$};
\node at (-.8,-1.7){\large$\beta$};
\node at (2,-1.5){\large$\alpha$};
\node at (1.4, 1.){\large$\gamma$};
\etp

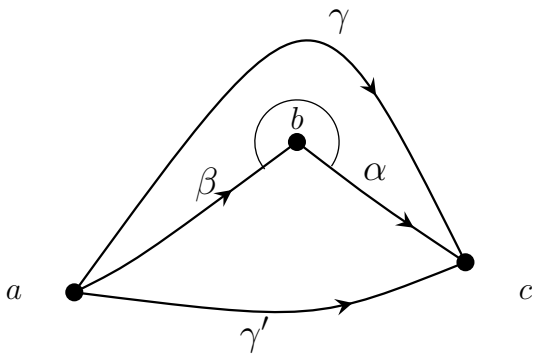
\caption{ The Picard-Lefschetz situation.  }
\label{fig:PicLef}
\enf
In this case we have identifications
\be\label{eq:3-isotop}
\Phi_{a,\gamma} \to \Phi_{a,\beta} \to \Phi_{a,  \gamma'}, \quad
\Phi_{c,\gamma'}\to \Phi_{c, \alpha}\to  \Phi_{c,\gamma},
\quad \Phi_{b, \beta} \to \Phi_{b, \alpha},
\ee
given by {\em clockwise} monodromies of the local systems $\bPhi$ around the corresponding arcs in the circles of directions.
So after these identications we can assume that we deal with single vector
spaces denoted by $\Phi_1, \Phi_3$ and $\Phi_2$ respectively.
Then we have \cite[Prop.1.8]{KS}  \cite[Prop.1.1.12]{KSS}:

\begin{prop}[Abstract Picard-Lefschetz identity]\label{prop:PL-form}
We have the equality of linear operators $\Phi_1\to\Phi_2$:
\[
m_{ac}(\gamma') = m_{ac}(\gamma) - m_{bc}(\alpha) m_{ab}(\beta).
\] \qed
\end{prop}
 

  \subsection {(Localized) perverse  sheaves on $\CC$}\label{subsec:locper}
  
  \paragraph{The category $\ol\Perv(X,A)$.} 
  
  We start with the case of an arbitrary Riemann surface, i.e. let 
  $X,A$ be as before. 
  Let $\LS(X)$ be the category of local systems
  of $\k$-vector spaces on $X$.  For each $\Lc\in\LS(X)$ the shifted sheaf
  $\Lc[1]$ is an object of $\Perv(X,A)$. It is straightforward that this defines
  an embedding of the shifted category $\LS(X)[1]$\
  (identified with $\LS(X)$) as
 a Serre subcategory on $\Perv(X,A)$ and so we have  the quotient abelian category
  \[
  \ol\Perv(X,A) = \Perv(X,A)/(\LS(X)[1]).
  \]
    Explicitly, $\Ob \ol\Perv(X,A)=\Ob\Perv(X,A)$
  while 
  \[
  \Hom_{\ol\Perv(X,A)}(\Fc, \Gc) =\Hom_{\Perv(X,A)}(\Fc, \Gc)/I_{\Fc, \Gc},
  \]
  where $I_{\Fc, \Gc}$ is the subset (actually a $\k$-vector subspace)
  formed by morphisms factoring as $\Fc \to\Lc[1]\to\Gc$ for some $\Lc\in\LS(X)$. 
  Thus we have a functor
   $\Perv(X,A)\to\ol\Perv(X,A)$ bijective on objects
  and surjective on morphisms. 
  
  Note that the first (but not the second) functor in \eqref{eq:PhiPsi-LS} 
  vanishes on $\LS(X)[1]$ and so  descends to a functor which we denote
  by the same symbol: 
  \be\label{eq:Phi-on-loc}
  \bPhi_a: \ol\Perv(X,A) \lra \LS(S^1_a). 
  \ee
  Further, let $\gamma$ be a simple path joining $a,b\in A$ as
  in \S \ref{par:transport}. The transport map $m_{ab}(\alpha)$ from
  \eqref{eq:transport} descends to  a natural transformations (also denoted
  $m_{ab}(\alpha)$)  between 
  functors \eqref{eq:Phi-on-loc} evaluated on $\dirr_a(\alpha)$ and $\dirr_b(\alpha)$.
  These transformation satisfy deformation invariance 
  \eqref{eq:def-inv} and the Picard-Lefschetz identities (Proposition \ref{prop:PL-form}).

   \paragraph{The Gelfand-MacPherson-Vilonen  description 
   of $\Perv(\CC,A)$ and $\ol\Perv(\CC,A)$ for finite $A$.} \label{par:GMV}
   
      From now on we assume that our Riemann surface $X$ is the complex line $\CC$
   with coordinate $w$.  We further assume that the set $A$ is finite. 
   
   \vskip .2cm
   
   Recall the descriptions  given in 
  \cite{gelfand-MV}.
From the topological point of view, we can replace $\CC$ by an open disk $D$ which 
we view as the interior of a closed disk $\ol D$, i.e.,  $D =  \ol D\-\del \ol D$, so 
$A=\{a_1.\cdots, a_N\} \subset D$. 
 
\begin{figure}[h]
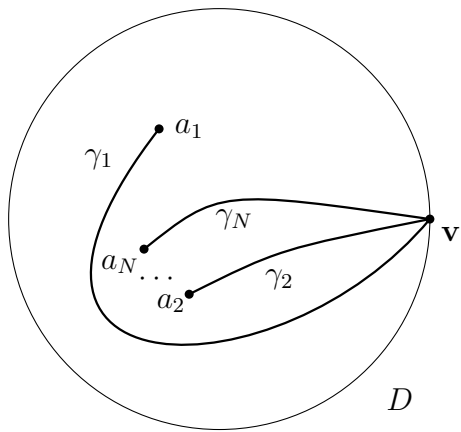

\centering
\btp[scale=.4, baseline=(current  bounding  box.center)]

\node (1)  at (-1, -2.5){};
\fill (1) circle (0.15);

\node (2)  at (-2.5, -1){};
\fill (2) circle (0.15);

\node (n)  at (-2,3){};
\fill (n) circle (0.15);

\draw (0,0) circle (7cm);

\node (b) at (7,0){};
\fill (b) circle (0.15);

\draw[line width = .3mm] (7,0) .. controls (2, -1) .. (-1, -2.5);

\draw[line width = .3mm] (7,0) .. controls (0,1) .. (-2.5,-1);

\draw[line width = .3mm] (7,0) .. controls (2,-6) and (-9,-6)   .. (-2,3);

\node at (7.7, -.5) {$\bv$};

\node at (6,-6) {$D$};

\node at (-1.6, -2.8) {$a_2$};

\node at (-2,-2){$\cdots$};

\node at (-3.3, -1.5){$a_N$};

\node at (-1, +3){$a_1$};

\node at (-4,+2){$\gamma_1$};

\node at (2,-2){$\gamma_2$};

\node at (0.5,0){$\gamma_N$};

\etp

\caption{  A spider defining the GMV-equivalence.  }
\label{fig:disk-2}
\end{figure}

  Fix a point $\bv\in \del D$. 
  Call a   $\bv$-{\em spider} for $(D,A)$   a system
$K = \{\gamma_1,\cdots, \gamma_N\}$ of simple closed  piecewise smoorth arcs in $D$ so that 
$\gamma_i$ joins $\bv$ with $a_i$ and
different $\gamma_i$ do not meet except at $\bv$, see   Fig. \ref{fig:disk-2}.

A $\bv$-spider for $(D,A)$ defines a total order on $A$, by  clockwise
ordering of the slopes of the $\gamma_i$ at $\bv$ (assumed distinct) and we choose
the numbering
$A=\{a_1, \cdots, a_N\}$  in this order, as in Fig. \ref{fig:disk-2}.

 Denote by $\Qc_N$  the category  of diagrams  of finite-dimensional $\k$-vector spaces
\[
\xymatrix{\Phi_N
\ar@<.4ex>[rd]^{u_N}
&
\\
\vdots & \Psi \ar@<.4ex>[ul]^{v_N}
\ar@<.4ex>[dl]^{v_1}
\\
\Phi_1 \ar@<.4ex>
[ur]^{u_1}
}
\]
such that   $T_{\Phi_i}:= \Id_{\Phi_i}-v_iu_i$ is invertible for each $i$. This implies that
each $T_{i, \Psi} = \Id_\Psi-u_iv_i$ is  invertible as well.

Given a  $\bv$-spider $K$  for $(D,A)$ and $\Fc\in\Perv(D,A)$, we consider, for 
each  $i=1,\cdots, N$, the space $\Phi_{i, K}(\Fc)= \Phi_{a_i,\gamma_i}(\Fc)$, i.e., 
the stalk of the local system $\bPhi_{a_i}(\Fc)$ at the point in
 $S^1_{a_i}$ represented by the direction of
$\gamma_i$.
Let us also identify the stalk $\Psi_{\gamma_i}(\Fc)$ with the stalk $\Fc_{\bf v}$,
 using the monodromy along $\gamma_i$. 
This gives a diagram
\[
\Theta_K(\Fc) \,=\, ( \Psi(\Fc),  \Phi_{i, K}(\Fc),  u_{i, \gamma_i} , v_{i, \gamma_i}  ) \,\in \, \Qc_N.
\]
Here $u_{a_i, \gamma_i}$ and $v_{a_i, \gamma_i}$ are the canonical maps along
 $\gamma_i$ as in
\eqref{eq:transport}.

\begin{prop} [\cite {gelfand-MV}] \label{prop:GMV1}
Let $K$ be a $\bv$-spider for $(D,A)$.
The functor $\Theta_K: \Perv(D,A)\to\Qc_N$ is an equivalence.
\qed
\end{prop}

Further,  a spider $K$ defines, for each $i\neq j$, a path $\alpha^K_{ij}$
joining $a_i$ and $a_j$ by first going from $a_i$ to $\bv$ and then from
$\bv$ to  $a_j$.

 Let $\Mc_N$ be the category whose objects are   diagrams consisting of:
\begin{itemize}
\item [(0)] Vector spaces $\Phi_i$, $i=1,\cdots, N$.

\item[(1)] Linear operators $m_{ij}: \Phi_i\to \Phi_j$ given for all $i,j$ (including $i=j$)
  such that $ \Id_{\Phi_i}-m_{ii}$ is invertible.
\end{itemize}

Given an object $\Fc\in\ol\Perv(D,A)$ and a spider $K$ as above,  we construct a diagram 
$\Xi_K(\Fc) = (\Phi_i, m_{ij}) \in\Mc_N$  by putting
$\Phi_i = \Phi_{i, \gamma_i}(\Fc)$ as before and
\[
m_{ij} \,=\begin{cases}
m_{ij}(\alpha^K_{ij}), &\text { if } i\neq j;
\\
\Id-T_i(\Fc), & \text { if } i=j,
\end{cases}
\]
 
\begin{prop}[\cite{gelfand-MV}]\label{prop:GMV2}
The functor $\Xi_K: \ol\Perv(D,A)\to \Mc_N$ is an equivalence. \qed

\end{prop}

\begin{rem}\label{rem:vladi}
 We can think of the point $\bv\in\del D$ as being  far away at the infinity (``Vladivostok'').
 This becomes even more natural if we use $D$ as a model for the compelx plane $\CC$. 
For this reason we will sometimes refer to the description of $\ol \Perv(D,A)$ given by Proposition
\ref{prop:GMV2} as the {\em Vladivostok description} and call the path
$m_{ij}^K$ the {\em Vladivostok path} joining $a_i$ and $a_j$. From the
naive ``physical'' point of view this is not the most natural way to connect $a_i$ and $a_j$
by a path. 

\end{rem}

   \paragraph{$\ol\Perv(\CC,A)$ inside $\Perv(\CC,A)$.} For
   $\Fc\in\Perv(\CC,A)$ the hypercohomology $H^i(\CC,\Fc)$
   of $\CC$ with coefficients in $\Fc$ vanish for $i\neq -1,0$. 
   Let $\Perv^0(\CC,A)\subset\Perv(\CC,A)$ be the full subcategory formed by
   $\Fc$ such that $H^i(\CC,\Fc)=0$ for all $i$.    The following  statement is 
   a reformulation of the  results of \cite{gelfand-MV, KKP}.

   \begin{prop}
   (a) The localization functor $\Perv(\CC,A) \to\ol\Perv(\CC,A)$ restricts to an equivalence
   of categories $\Perv^0(\CC,A)\to\ol\Perv(\CC,A)$. 
   
   \vskip .2cm
   
   (b) Each object of $\Perv^0(\CC,A)$ reduces to a single sheaf in degree $(-1)$. 
   \end{prop}
   
   \noindent{\sl Proof:}  (a) In the proof of \cite[Prop. 2.3]{gelfand-MV}  the authors
   construct   a full embedding   $\lambda: \Mc_N \to \Qc_N$,
  where $\Mc_N$ is the category of diagrams describing
  $\ol\Perv(\CC,A)$ by Proposition \ref{prop:GMV2} and $\Qc_N$ is the
  category of diagrams describing $\Perv(\CC,A)$
  by Proposition \ref{prop:GMV1}.  The image $\Im(\lambda)$ is a subcategory
  in $\Qc_N$ which maps equivalently to $\Mc_N$ under the localization
  functor $\Qc_N\to\Mc_N$.  It is then verified directly from the definition
  of $\lambda$ , see the end of the proof of \cite[ Thm.2.29]{KKP} that $\Im(\lambda) = \Perv^0(\CC,A)$. 
  
  \vskip .2cm
  
  (b)  This is shown in the first part of the proof of 
   \cite[ Thm.2.29]{KKP} . \qed

   \begin{rem}
   Se we can write the lifting functor $\lambda$ in a geometric way, as
    \[
   \lambda: \ol\Perv(\CC,A)\buildrel \sim\over  \lra \Perv^0(\CC,A)\subset\Perv(\CC,A). 
   \]
   By Proposition  \ref{prop:GMV2}, an object of
   $\ol\Perv(\CC,A)$ is determined by its vanishing cycles $\Phi_a$  and
   transport maps between them, but does not have  well defined  stalks 
   at points outside
    $A$, which for an actual perverse sheaf form
   a (shifted)  local system on $\CC\-A$. The functor $\lambda$ provides
   a preferred way to supply such local system. 
   Interpreting the construction of $\lambda$ from
   \cite{gelfand-MV} in a geometric fashion, we see that 
    the typical stalk of this system is identified
   with $\bigoplus_{a\in A} \Phi_a$. 
   \end{rem}
   
   \begin{ex}\label{ex:unit-obj}
   let $a\in A$  and $\Fc = \k_a\in\Perv(\CC, A)
   $ be the skyscraper
   sheaf at $a$. Denoting $\ol\Fc\in\ol\Perv(\CC,A)$ the image of $\Fc$,
   the lift $\lambda(\ol\Fc)\in\Perv^0(\CC,A)$ is the sheaf
   $j_!( \ul\k_{\CC\-\{a\}})[1]$, where $j: \CC\-\{a\}\to\CC$ is the embedding. 
   It becomes isomorphic to $\Fc=\k_a$ in $\ol\Perv(\CC, A)$ because
   of the exact sequence of sheaves
   \[
   0\to j_!( \ul\k_{\CC\-\{a\}})  \to \ul\k_\CC \to  \k_a \to 0
   \]
   gives an exact sequence in $\Perv(\CC,A)$
   \[
   0\to \k_a \to  j_!( \ul\k_{\CC\-\{a\}})[1]  \to \ul\k_\CC [1]\to 0
   \]
   with third term in $\LS(\CC)[1]$. 
   \end{ex}
   
   
   \subsection{Fourier transform of perverse sheaves 
   and their Stokes data} \label{subsec:FT}
   
   In this section we assume $\k=\CC$. 
   
   \paragraph {  The formal Fourier transform.} 
  Let  $D_w = \CC\langle w, \del_w\rangle$ be the Weyl algebra of polynomial differential operators on
  $\CC$ and ${D_w}\Mod^\hol\supset D_w\Mod^{hrs}$ be the categories of holonomic and
  holonomic regular singular $D_w$-modules.  It is well known that the {\em solution functor}  (on all holonomic modules, regular singular or not)
     \[
   \Sol: D_w\Mod^\hol \lra \Perv(\CC), \quad M\mapsto \Sol(M) = \ul\RHom_{D_w}(M, \Oc_\CC)[1]
   \]
   takes values in $\Perv(\CC)$. Further,  its restriction to $D_w\Mod^{hrs}$ 
   is an equivalence (Riemann-Hilbert correspondence). See, e.g.,  \cite{malgrange}. 
    So we can realize any $\Fc\in \Perv(\CC)$ as $\Sol(M)$ for a unique
   $M = M_\Fc\in D_w\Mod^\hr$.  
      
  The {\em formal Fourier transform} is 
 the isomorphism
 \[
 D_w = \CC\langle w, \del_w\rangle \lra D_z = \CC\langle z, \del_z\rangle, \quad w\mapsto -\del_z, \, \del_w\mapsto  z,
 \]
 matching the analytic Fourier transform of solutions. Given $M\in D_w\Mod^\hol$, its {\em Fourier transform} 
 $\wh M$ is the same $M$ but considered as a $D_z$-module using the above isomorphism.
 We refer to 
  \cite{malgrange} for background on this construction. In particular, it is known 
   that $ \wh M$ is again holonomic, 
   so  we have the perverse sheaf $\Sol(\wh M) =  \ul\RHom_{D_z}(M, \Oc_\CC)[1]$. 
   
  If $M$ is regular, then  $\wh M$ is typically not regular. In this case it  is also known that 
   the  perverse sheaf $\wh\Fc = \Sol(\wh M)$
 has $0$ as  the only possible  singularity,  so we get the  functor
 (also called the {\em Fourier transform}) 
 \be
 \FT: \Perv(\CC) \lra \Perv(\CC,0), \quad \Fc\mapsto  \wh \Fc := \Sol(\wh M_\Fc). 
 \ee
 The restriction of $\FT(\Fc)$  to $\CC\-\{0\}$  has thus the form $\FT_\gen(\Fc)[1]$
 for a local system $\FT_\gen(\Fc)$ on $\CC^*$, so we have the functor
 \[
 \FT_\gen: \Perv(\CC) \lra \LS(\CC^*). 
 \]
 
   As $\CC^*$ is homotopy equivalent to the unit circle $S^1\subset\CC$,
 for any local system $\Lc$ on $\CC^*$ we can speak about the stalk $\Lc_\zeta$
 at ant $\zeta\in S^1$. Further,  for any $a\in \CC$, the circle $S^1_a$ 
 of directions at $a$ is identified with this fixed $S^1$, so we
 can view any $\bPhi_a(\Fc)$ as a local system on this $S^1$.
 With this understanding, the following is true \cite[Ch.XII]{malgrange}
 \cite [Prop.6.1.4]{dagnolo-sabbah}.

 \begin{prop}\label{prop:FTgen}
 Let $\Fc\in\Perv(\CC)$. We have a natural isomorphism
\[
\FT_\gen(\Fc)\,  \= \, 
 \bigoplus_{a\in\CC} \bPhi_a(\Fc)
 \]
  of local systems on $S^1$.    \qed
 \end{prop}

 In particular, the functor $\FT_\gen$ factors through $\ol\Perv(\CC)$,
 which is one of the reasons to consider this localized category. 
 
 \vskip .2cm
 
\begin{figure}[h]
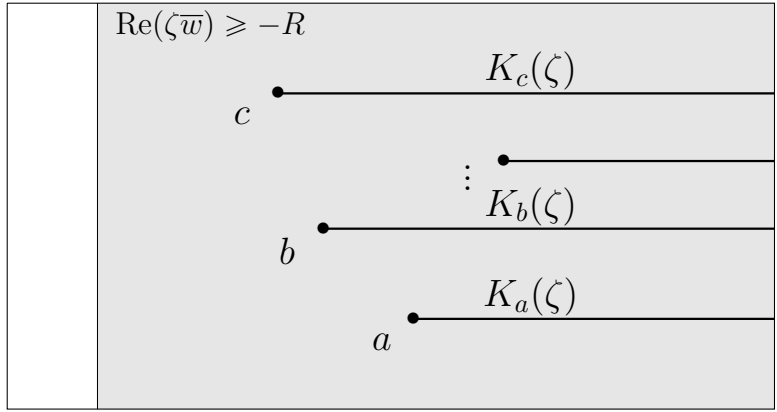

\centering
\btp[scale=.6]

\draw (-8, 4) -- (-8,-5) -- (9,-5) -- (9, 4) -- (-8,4);
\node at (1,-3){$\bullet$};
\node at (-1,-1) {$\bullet$};
\node at (3,0.5) {$\bullet$};
\node at (-2,2) {$\bullet$};


\draw [line width=.3mm] (1,-3) -- (9,-3);
\draw [line width=.3mm] (-1,-1) -- (9,-1);
\draw [line width=.3mm] (3,0.5) -- (9,0.5);
\draw [line width=.3mm] (-2,2) -- (9,2);

\node at (0.3, -3.5){\large$a$};  \node at (3.6,-2.5){\large$K_a(\zeta)$};
\node at (-1.8,  -1.5){\large$b$};  \node at (3.6, -0.5){\large$K_b(\zeta)$};
\node at (2.2,  0.3){\large$\vdots$};
\node at (-2.8,  1.5){\large$c$};    \node at (3.6, 2.5){\large$K_c(\zeta)$};

\filldraw[opacity=0.1] (-6,4) -- (-6,-5) -- (9,-5) -- (9, 4) -- (-6, 4);
\draw[line width= 0.1mm] (-6, 4) --(-6,-5);

\node at (-3.5, 3.5){$\Re(\zeta \ol w)\geq -R$};
 
\etp

\caption{The half-lines $K_a(\zeta)$.  Here $\zeta=1$. } \label{fig:horiz-cuts}
\end{figure}

 \noindent{\sl Proof:} For future reference, we give a sketch of the construction
 of the identification of the proposition for generic $\zeta$. 
 Let $\Fc\in\Perv(\CC, A)$. A general result  \cite[Th.3.1.1]{daia}
 expresses the stalk of $\FT(\Fc)$  at $\zeta\in S^1$ (considered as the unit circle in $\CC$)
 as the cohomology with support,namely
 \be\label{eq:FTFC-zeta}
 \FT(\Fc)_\zeta = H^0_{\{ \Re(\zeta \ol w)\geq -R\}} (\CC, \Fc), \quad R\gg 0, 
  \ee
 with $R$ large enough so that the shifted half-plane contains $A$. For $a\in A$ let
  $K_a(\zeta) = a + \zeta\cdot \RR_+$ be the half-line in the direction $\zeta$ issuing from
  $a$,   see Fig. \ref{fig:horiz-cuts}.

  The embedding of the complements
  \[
  \CC \-  \,\, \{ \Re(\zeta \ol w)\geq -R\} \,\,\hra \,\, \CC \- \,\,\, \bigcup _{a\in A} K_a(\zeta)
  \]
  is a homotopy equivalence and $\Fc$ is locally constant on both of them. 
 So  we can  replace the support in \eqref{eq:FTFC-zeta} by the union  $\bigcup _{a\in A} K_a(\zeta)$. 
   Further, if $\zeta$ is generic enough, then the union is disjoint
   and so we have
  \[
  \FT(\Fc)_\zeta \,\= \, H^0_{\bigsqcup_{a\in A}  K_a(\zeta)}(\CC, \Fc) \, =\, \bigoplus_{a\in A} H^0_{K_a(\zeta)}
  (\CC, \Fc) \,=\, \bigoplus_{a\in A} \Phi_{a} (\Fc). \qed
  \]

 \paragraph {The Stokes filtration of $\wh\Fc$.}
 As $\wh M$ is typically irregular, its solutions (i.e., sections of $\FT_\gen(\Fc)$)
 can grow exponentially in sectors near $\oo$. So
 the local system $\FT_\gen(\Fc)$ carries
 the additional {\em Stokes structure} given by the data  of such exponential growth
 in  various  sectors.  In our case, the growth is at most of the form $e^{\lambda z}$
 where $\lambda$ is a constant. 
 
 Thus  each stalk $\FT(\Fc)_\zeta)$ carries the {\em Stokes filtration} 
 $(\Sigma _\lambda)_{\lambda\in \CC}$
labelled by   the set $\CC$ with partial order $\leq_\zeta$ given by $\lambda \leq_\zeta \mu$
 if $\Re(\zeta \lambda) \leq \Re(\zeta \mu)$, i.e., $e^{\lambda z}$ is dominated  by $e^{\mu z}$, as $z\to\oo$
 on the ray $\zeta\cdot \RR_+$. The subspace $\Sigma_\lambda \FT(\Fc)_\zeta
 \subset \FT(\Fc)_\zeta$ consists of solution of $\wh M$ which grow in the direction
 $\zeta$ at most as $e^{\lambda z}$. 
 
 \begin{prop}\label{prop:stokes-filt}
 For a generic $\zeta$ we have
 \[
 \Sigma_\lambda \FT(\Fc)_\zeta \,=\, \bigoplus_{\Re(\zeta a) \geq\ -\lambda} \Phi_a(\Fc). 
 \]
 \end{prop}
 
 \noindent{\sl Proof: }This follows from representing solutions of $\wh M$ as actual
 Fourier integrals $g_{i, \zeta}(\phi)(z) $ corresponding to $\phi\in \Phi_{a_i}(\Fc)$
  over the half-lines $K_i(\zeta)$ as in Fig. \ref{fig:horiz-cuts}.
  See \cite[Ch. XII]{malgrange} for more details. 
 The  growth near $\zeta \oo$  of  $g_{i, \zeta}(\phi)(z)$ is of the rate  $e^{-a_i z}$. 
 \qed
 
 
 \subsection{The Lefschetz perverse sheaf and its Fourier transform} \label{subsec:lefper}
 
 In this section we take $\k$ to be an arbitrary field.

 \paragraph{The exponential integral as Fourier transform.}  Let $X$ be a smooth
  complex algebraic variety of dimension $n$   and
 $S: X\to\CC$ be a regular  function.  Given a regular volume form $dx$ on $X$, we can
 consider exponential integrals
 \be\label{eq:exp-int}
 I(z) = I(\hbar) = \int_{\Gamma} e^{{i\over\hbar} S(x)} dx, \quad z=1/\hbar, 
 \ee
 where $\Gamma$ is a locally finite   $n$-cycle in $X$ such that the integral converges.  
 We can consider $I(z)$ as a multivalued function whose determinations are labelled
 by choices of $\Gamma$. 
 
 It is classical that one can split the integration in  $I(z)$ into two stages. First, we have the relative volume
 form $dx/dS$, a rational section of $\Omega^{n-1}_{X/\CC}$ and form the multivalued
 function $L_S$ on $\CC$ 
 \[
 L_S(w) = \int_{\Gamma_w\subset S^{-1}(w)} {dx\over dS}
 \]
 where $\Gamma_w$ is an $(n-1)$-cycle in the fiber varying with $w$ via the Gauss-Manin connection. 
Then formally (assuming  that $\Gamma$ is formed out of the $\Gamma_w, w\in\gamma$ for a cycle
$\gamma$ in $\CC$)
\[
I(z) = \int_{ \gamma} L_S(w) e^{izw} dw
\]   
is (a determination of) the Fourier transform of $L_S$.  Note that $L_S$   is a Nielsen type function,
satisfying a differential equation with regular singularities, so all transcendental nontriviality of
$I(\hbar)$ comes from the $1$-dimensional Fourier transform. 

\paragraph{ The Lefschetz perverse sheaf.} \label{par:lefper}
A sheaf-theoretic analog of the function
$L_S$ is given by  the collection of perverse sheaves
\[
\Lc_S ^i= \ul H^i_{\perv} (RS_* (\ul \k_X[\dim(X)]), \quad i\in\ZZ
\]
on $\CC$. Here $\ul H^i_\perv$ is the degree $i$ perverse cohomology taken  with
respect to the perverse t-structure.  We will be particularly interested in the case $i=0$
and write simply $\Lc_S = \Lc^0_S$ while using the notation $\Lc_S^\bullet$ for
the graded perverse shead $\bigoplus_i \Lc_S^i$. 

\begin{prop}
Assume that the function $S: X\to \CC$ has only isolated singularities and is proper as a morphism
of algebraic varieties (so each $S^{-1}(a)$ is compact). Then: 
\vskip .2cm

(a) If $a\in\CC$ is a non-critical value for $S$, then $\Lc_S[-1]$ is a local system near $w$ and its
stalk at $a$ is identified with $H^n(S^{-1}(a), \k)$.

\vskip .2 cm

(b) If $a$ is a critical value   then
\[
\Phi_a(\Lc_S) \,=\, \bigoplus_{x\in S^{-1}(a)} \Phi_S(\ul \k_X[\dim(X)])_x
\]
is the direct sum of the classical (Lefschetz) 
spaces of vanishing cycles for $S$ at the critical points over $a$. 
\end{prop}

\noindent Thus, for $\k=\CC$ sections of the local system $\Lc_S[-1]$ over the open set
of non-critical values of $S$  give determinations of $L_S$. 

\vskip .2cm 

\noindent {\sl Proof:} (a) follows since the perverse t-structure is centered around the middle dimension. 
Part (b) follows from the definition of  vanishing cycles. \qed

\begin{rem}
One can say that the first step towards categorification of the exponential integral \eqref{eq:exp-int} is given
 by $\FT(\Lc_S)$, the Fourier
transform of the perverse sheaf $\Lc_S$. \footnote{For the next step one should replace cohomology groups by appropriate categories, cf. \cite{KSS}.}As we saw in \S \ref{subsec:FT}, the
structure of $\FT(\Lc_S)$ near $\oo$ is entirely given by 
the image of $\Lc_S$ in $\ol\Perv(\CC)$, i.e.,  by the vanishing cycles of $\Lc_S$
and the transport maps $m_{ab}(\gamma)$ between then. To find  these data, we do
not need to compactify $S$ to a proper morphism $X\to\CC$. 
For example, if $S$ is a Morse function, 
 then  it suffices to have
the part of $X$ containing the critical points and the  Lefschetz thimbles emanating from critical points towards other critical points. 
\end{rem}


 \section{ Alien derivatives for perverse sheaves: elementary theory}\label{sec:alien}

\subsection{Additive  convolution of localized perverse sheaves and Fourier transform}
\label{subsec:conv}

 Convolution of  \'etale perverse sheaves on commutative algebraic
  groups was studied by N. Katz \cite{katz}. We will need a simplified
  version for analytic perverse sheaves on the  additive group $\CC$,
see e.g. \cite {fressan-jossen}.  
   
  \paragraph{ Additive convolution. }  
  Let $\k$ be an aribtrary field.
    Let $A,B\subset\CC$ be finite subsets and
   let $\Fc\in\Perv(\CC,A)$, $\Gc\in\Perv(\CC,B)$. Then 
   $\Fc\boxtimes\Gc$ is a constructible complex (in fact, a perverse sheaf)
   on $\CC\times\CC$. 
   We have the addition map
   \[
   + : \CC \times\CC\lra \CC,  \quad (w', w'')\mapsto w'+w''. 
   \]
   Let $A+B= +(A\times B)$ be the set of sums $a+b$, $a\in A$, $b\in B$.
   The map $+$
    gives  the {\em additive convolution} 
      \[
     \Fc * \Gc := R(+)_*(\Fc\boxtimes\Gc) 
   \]
   Note that $ \Fc *\Gc$  is  a priori a constructible complex
   on $\CC$ with singularities (points of local non-constancy of cohomology)
   contained in $A+B$. 
   
   \begin{prop}\label{prop:F*G}   \cite[Prop.2.4.3] {fressan-jossen}
    If $\Fc\in \Perv(\CC,A)$ and $\Gc\in\Perv^0(\CC,B)$, then
   $\Fc* \Gc\in \Perv^0(\CC, A+B)$. \qed
    \end{prop}
   
   Let $\Perv(\CC) = \bigcup_A \Perv(\CC,A)$ be the category of all perverse sheaves on $\CC$
   with finitely many singularities. In a similar way we define the quotient category
    $\ol\Perv(\CC)$ and  its lift  $\Perv^0(\CC)\subset\Perv(\CC)$. By the above, $\ol\Perv(\CC)$
    and $\Perv^0(\CC)$ are equivalent.  
    
    \begin{cor} \cite[Th.2.4.11]{fressan-jossen}
    The operation $*$ makes $\ol\Perv(\CC)\simeq  \Perv^0(\CC)$ into a symmetric monoidal
    category with unit object $\1$ being the class of  $\ul\k_0$ in $\ol\Perv(\CC)$ or, equivalently,
    its lift $j_! \ul\k_{\CC\-\{0\}} \in \Perv^0(\CC)$, see Example \ref{ex:unit-obj}. 
    \end{cor}
   
   \paragraph{Comparison with the  Hurwitz convolution for analytic functions.} 
    
  The operation $*$ for perverse sheaves is a categorical analog of the additive (or Hurwitz)
   convolution of holomorphic functions $f,g\in\CC\{\{w\}\}$ defined near $0$:
   \be\label{eq:hurwitz}
   (f*g)(w) = \int_0^w f(u) g(w-u) du. 
   \ee
   If
    \[
   f(w) = \sum_{n=0}^\oo a_n {w^n\over n^!}, \quad g(w) = \sum_{n=0}^\oo b_n{w^n\over n!},
   \]
   then 
   \[
   (f*g)(w) = \sum_{n=0}^\oo \biggl(\sum_{i+j=n} a_i b_j\biggr) {w^{n+1}\over (n+1)!}
   \]
   (series without constant term, in fact we have $1*1=w$).  The theorem of Hurwitz
   (additive version of the Hadamard theorem of multiplication of singularities)
 says that if $f,g$ extend to possibly multivalued analytic
 functions in $\CC$ with singularities in possible infinite sets $A,B$ respectively,
 then $f*g$ similarly extends to a possibly multivalued analytic function in $\CC$
 with singularities in $A\cup B\cup (A+B)$. See \cite{hurwitz}  and
 later treatments in \cite {schottlaender} and \cite{sauzin} \S6.4. 
 The additional possible singularities at $A = A+\{ 0\}$ and $B=\{ 0\}+B$ in Hurwitz's
 theorem as compared to Proposition \ref{prop:F*G} come from the 
 fact that the integration path in \eqref{eq:hurwitz}, starting from $0$,
 is a chain but not a cycle with coefficients in the local system of determinations of the integrand.

 \paragraph{Additive convolution and Fourier transform.}

 Let $\k=\CC$. 
  The classical principle  that ``Fourier transform takes convolution into product''
  has  in 
 our case  the following form.

 \begin{prop}\label{prop:FTvs*}
 For $\Fc, \Gc\in \Perv^0(\CC)$ we have a natural isomorphism of local systems on
 $S^1$
 \[
 \FT_\gen(\Fc * \Gc) \, \simeq \, \FT_\gen(\Fc) \otimes\FT_\gen(\Gc). 
 \] 
 \end{prop}
 
\noindent{\sl Proof:} By the Riemann-Hilbert correspondence, the 
 direct image $R+_*$ in the definition of $\Fc * \Gc$ can be calculated at the level
 of $D$-modules. That is, let $M,N\in D_w\Mod^{hrs}$. We define the
 $D$-module additive convolution to be the complex of $D_w$-modules
 \[
  M*^D N = R+^D_*(M\boxtimes N),
  \]
  where:
  \begin{itemize}
  
  \item [(1)]  $M\boxtimes N = M\otimes_\CC N$ considered as a module
 over $D_w\otimes_\CC D_w = \CC\langle w', w'', \del_{w'}, \del_{w''}\rangle$
 (polynomial differential operators on $\CC\times \CC$);
 
 \item [(2)]  $R+_*^D$ is the derived $D$-module direct image, i.e., the de Rham complex
 along the fibers of $+:\CC\times \CC\to \CC$. 
  
 \end{itemize}  
 
 So defined, $M*^D N$  is a complex with holonomic regular singular
 cohomology  and $\Sol(M*^D N) \simeq  \Sol(M) * \Sol(N)$
 in the derived category of  constructible complexes.

 Further, the tensor product of local systems in the claim of the proposition
 corresponds to the tensor product of $D_z$-modules over $\CC[z]$. So our claim
 is a consequence of the following one.
 
 \begin{lem}
 Let $M,N\in D_w\Mod^{hrs}$. Then we have a natural isomorphism 
 in derived category
 \[
 \wh{M*^D N} \= \wh M \otimes^L_{\CC[z]} \wh N
 \]
 \end{lem}
 
 \noindent{\sl Proof of lemma:}  The fibers of the map $+$ are $1$-dimensional
 with the relative tangent bundle generated by the vector field
 $\del_{w'}-\del_{w''}$. So 
 $M*^D N$ is the complex
 \[
 M\otimes_\CC N \buildrel \del_{w'}-\del_{w''}\over \lra M\otimes_\CC N
 \]
 with the differential $\del_{w'}-\del_{w''} = \del_w\otimes 1 - 1\otimes\del_w$.
 Now,   $\wh M \otimes^L_{\CC[z]} \wh N$ is the complex
 \[
 \wh M\otimes_\CC \wh N \buildrel z\otimes 1 - 1\otimes z\over\lra \wh M \otimes_\CC \wh N.
 \]
But  $\wh M$ is $M$ in which $z$ acts as $\del_w$ and $\del_z$ as
 $-w$, and similarly for $\wh N$, so the second complex is identified
 with the first after the  Fourier transform. 
 This proves the lemma and Proposition \ref{prop:FTvs*}.
 \qed

 \paragraph{The Thom-Sebastiani theorem.}  Propositions  \ref{prop:FTgen}
 and \ref{prop:FTvs*} lead to a multiplicativity property  (a version of the 
 Thom-Sebastiani theorem)
 which does not involve Fourier transform and can be proved directly
 at the level of perverse sheaves.  We take $\k$ to be an arbitrary field. 
 
 \vskip .2cm

 Let $\LS(S^1)^\CC$ be the category of {\em $\CC$-graded local systems} 
 on $S^1$, i.e. of collections $\Lc = (\Lc_c)_{c\in\CC}$ of local systems on $S^1$
 such that $\Lc_a=0$ for almost all $a$. 
 This category has a symmetric monoidal structure given by
 \[
 (\Lc\otimes\Mc)_c \, = \, \bigoplus_{a + b=c} \Lc_a \otimes\Mc_b. 
 \]
 We have the {\em total vanishing cycle functor}
 \[
 \bPhi: \Perv(\CC) \lra \LS^\CC(S^1), \quad \Fc \mapsto \bPhi(\Fc) :=  
 (\bPhi_c(\Fc))_{c\in \CC}. 
 \]

 \begin{thm} \label{thm:Thom}  \cite[Th.2.8.3] {fressan-jossen}
  The functor 
  $\bPhi$ is  symmetric monoidal. In other words, 
 for any $\Fc, \Gc\in\Perv^0(\CC)$ and  any $c\in \CC$
  we have a natural isomorphism of
 local systems on $S^1$
 \[
 \bPhi_c (\Fc * \Gc) \, \= \, \bigoplus_{a+b=c} \bPhi_a(\Fc) \otimes \bPhi_b(\Gc).  
  \]
 \end{thm}  
 
 For convenience of the reader and future reference we give a direct proof.   
 It is enough to consider the case $c=0$, the case of arbitrary $c$ is similar. 
 We further identify the stalks of the local systems of vanishing cycles
  at the point $1\in S^1$,
 i.e., in the direction of $\RR_+$, the identification of monodromy
 following by the same arguments as below.  We will use the letter $\Phi$ to mean such
 stalks.
 
\begin{figure}[h]
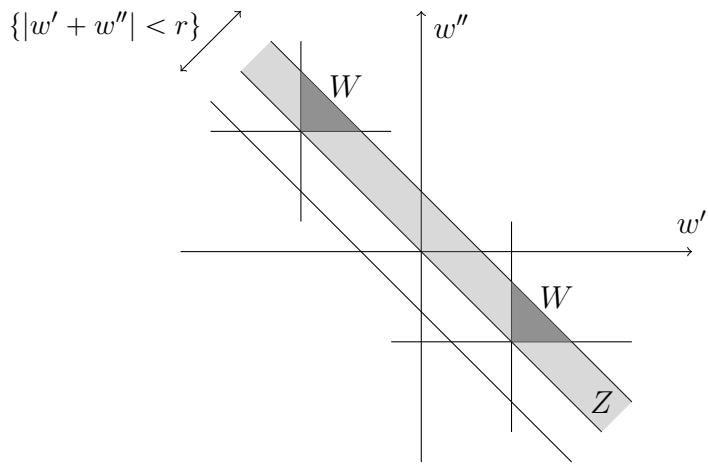

\centering
\btp[scale=.4, baseline=(current  bounding  box.center)]

\draw [->] (-8,0) -- (9,0); 
\draw [->] (0, -7) -- (0,8); 
\draw (-6,6) -- (6, -6); 

\draw (-7,5) -- (5, -7); 
\draw (-5,7) -- (7, -5); 

\draw (-7,4) -- (-1, 4); 
\draw (-4,1) -- (-4,7); 

\draw (-1, -3) -- (7, -3); 
\draw (3, -6) -- (3,1);

\draw [<->] (-8, 6) -- (-6,8); 

\node at (-10.5, 7.5) {\small $\{|w'+w''|<r\}$}; 

\fill[gray, opacity=0.3] (6, -6) -- (7, -5) -- (-5, 7) -- (-6,6) -- (6, -6) ; 

\node at (6, -5) {$Z$}; 

\fill[gray, opacity=0.8] (-4,4) -- (-2,4) -- (-4,6) -- (-4,4); 

\fill[gray, opacity = 0.8] (3, -3) -- (5, -3) --( 3, -1) -- (3, -3); 

\node at (-2.5, 5.5) {$W$}; 
\node at (4.5, -1.5) {$W$}; 
\node at (9,1) {$w'$}; 
\node at (1,7.5) {$w''$}; 
\etp

\caption{The two sets of supports in the Thom-Sebastiani theorem. }
\label{fig:thom}
\end{figure}

 For a perverse sheaf $\Ec\in\Perv(\CC)$ we have a canonical identification
 \[
 \Phi_0(\Ec) \,=\, R\Gamma_{\{\Re(w)\geq 0\}} (\{|w|<r\}, \Ec),
  \]
  where $r>0$ is small enough. Indeed, the discussion in \S \ref{subsec:perRiem} \ref{par:nearvan}
gives $\Phi_0(\Ec)$ as  $R\Gamma_{\RR_+}(\{|w|<r\}, \Ec)$ 
but the half-plane $\Re(w)\geq 0$ (the part of it lying in the disk $|w|<r$) 
is stratified homotopy equivalent (w.r.t. to the stratification of $|w|<r$ by $0$
and everything else) to the real half-line. 
Therefore, using $w', w''$ as coordinates in $\CC\times\CC$, we have
\be\label{eq:phi-FG}
\Phi_0(\Fc * \Gc) \, =\, R\Gamma_Z
(\{ |w'+w''| < r\}, \Fc\boxtimes\Gc), \quad Z:= \{ \Re(w' +w'') \geq 0\},
\ee
see Fig. \ref{fig:thom}. 
Now consider the subset 
\[
W \, := \, \bigcup_{a,b\in A, a+b=0} \bigl\{ \Re(w'-a)\geq 0, \, \Re(w'' -b)\geq 0\bigr\} 
\,\subset \, Z. 
\]

If  $r$ is small enough, then the union in the definition of $W$ is disjoint and
\be\label{eq:PhiF-PhiG}
R\Gamma_W(\{|w' +w''|< r\}, \Fc\boxtimes\Gc) \,=\, \bigoplus_{a,b\in A, a+b=0}
\Phi_a(\Fc)\otimes\Phi_b(\Gc)
\ee
by the K\"unneth theorem. It remains to notice that the complements of the two supports,
i.e., the open subsets
\[
U\,=\, \{|w'+w''|<r\}\,  \-\, Z \quad \text{and} \quad V \,=\, \{w' +w''|<r \} \,\-\, W
\]
are  stratified homotopy equivalent (with respect to the stratification given by
the singularities of $\Fc\boxtimes\Gc$). More precisely, 
$U$ can be obtained as a deformation retract of $V$ with deformations
affecting only the area where $\Fc\boxtimes\Gc$ is locally constant, so $R\Gamma_Z = R\Gamma_W$ and therefore
\eqref{eq:phi-FG} and \eqref {eq:PhiF-PhiG} give the same answer. \qed

\vfill\eject

 \subsection{Rectilinear transports with avoidances and alien transports}\label{subsec:trans-avoid}

\paragraph {Rectilinear transports.} \label{par:rect-trans}
We start by taking $\k$ to be an arbitrary field. 
Let $a,b\in A$ be two distinct points. The most natural  path
joining $a$ and $b$ is the straight line interval $[a,b]$.  We denote by
\be\label{eq:zeta-ab}
\zeta_{ab}\, = \, {b-a \over |b-a|}\, \in\,  S^1\, :=\,\{\zeta\in \CC: |\zeta|=1\}
\ee
 the slope of $[a,b]$ in the direction from $a$ to $b$.

If $[a,b]$ does not
contain any other elements of $A$, then for any $\Fc\in\Perv(\CC,A)$
we can use $[a,b]$ to define  the {\em rectilinear transport} 
\be
 m_{ab}= m^\Fc_{ab} = m^\Fc_{ab}([a,b]): \bPhi_{a, [a,b]}(\Fc) \lra \bPhi_{b, [a,b]}(\Fc). 
\ee
We say that $A$ is {\em in linearly general position including $\oo$}, if no three points of $A$
lie on a real line in $\CC$ and no $\zeta_{ab}$ belongs to  $\RR$. 

Suppose $A$ is in linearly general position including $\oo$. Then   all rectilinear transports
$m_{ab}([a,b])$, $a,b\in A$,  are defined. Let  $\Mc_N$  be the category of diagrams $(\Phi_i, m_{ij})$
as in Proposition 
  \ref {prop:GMV2}.
Let us number  $A=\{a_1,.\cdots, a_N\}$ in
an arbitrary way, and denote $\zeta_{ij}=\zeta_{a_i, a_j}$, $i\neq j$ and $\bPhi_i(\Fc) = \bPhi_{a_i}(\Fc)$
(a local system on $S^1$). 
In this notation, we define a functor
\[
\Xi_\rect: \ol \Perv(\CC,A) \lra \Mc_N, \quad \Fc\mapsto (\Phi_i, m_{ij}), 
\]
as follows. We put
\[
\Phi_i = \bPhi_{a_i, a_i+\RR}(\Fc) = \bPhi_i(\Fc)_1
\]
(the stalk of the local system $\bPhi_i(\Fc)$ at $1\in S^1$). Further, 
$m_{ii}$ is defined as $\Id-T_i(\Fc)$ where $T_i(\Fc): \bPhi_i(\Fc)_1\to\bPhi_i(\Fc)_1$ is the counterclockwise
monodromy. For $i\neq j$ we define $m_{ij}$ as the composition
\be\label{eq:m_{ij}-rect-GMV}
\bPhi_i(\Fc)_1\buildrel T_1^{\zeta_{ij}}\over  \lra \bPhi_i(\Fc)_{\zeta_{ij}} =\bPhi_{a_i,[a_i, a_j]}(\Fc) \buildrel m_{a_i, a_j}^\Fc\over\lra
\bPhi_{a_j, [a_j, a_i]}(\Fc) = \bPhi_j(\Fc)_{\zeta_{ji}} \buildrel T_{\zeta_{ji}}^1 \over\lra  \bPhi_j(\Fc)_1,
\ee
where $T_1^{\zeta_{ij}}$  (resp. $T_{\zeta_{ij}}^1$) is the monodromy of the local system $\bPhi_i(\Fc)$ from $1$ to $\zeta_{ij}$ 
(resp. from $\zeta_{ji}$ to $1$) taken in the
counterclockwise direction, if $\Im(w_i)< \Im(w_j)$ and in the clockwise direction, if $\Im(w_i) > \Im(w_j)$.

  \begin{prop}\label{prop:Locperv-gen}
  If $A$ is in linearly general position including $\oo$, then the functor $\Xi_\rect$
  is an equivalence of categories. 
 \end{prop}
 
 \noindent{\sl Proof:} This statement, which is  \cite[Prop.2.1.7]{KSS},
 is deduced from Proposition   \ref {prop:GMV2} by deforming the set $A$ to the convex position. 
  \qed

\begin{figure}[h]
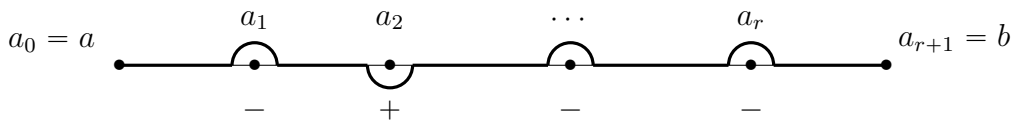

     \centering
\btp[scale=.6, baseline=(current  bounding  box.center)]

\draw (-8,0) -- (9,0); 

\node at (-8,0) {\small$\bullet$}; 
\node at (-5,0) {\small$\bullet$}; 
\node at (-2,0) {\small$\bullet$}; 
\node at (2,0) {\small$\bullet$}; 
\node at (6,0) {\small$\bullet$}; 
\node at (9,0) {\small$\bullet$}; 

\draw[line width=1.2] (-4.5,0) arc [start angle=0, end angle=180, x radius=0.5cm, y radius=0.5cm];
\draw [line width=1.2]  (-2.5,0) arc [start angle=180, end angle=360, x radius=0.5cm, y radius=0.5cm];
\draw [line width=1.2]  (2.5,0) arc [start angle=0, end angle=180, x radius=0.5cm, y radius=0.5cm];
\draw [line width=1.2]  (6.5,0) arc [start angle=0, end angle=180, x radius=0.5cm, y radius=0.5cm];

\draw[line width=1.2] (-8,0) -- (-5.5,0); 
\draw[line width=1.2] (-4.5,0) -- (-2.5,0); 
\draw[line width=1.2] (-1.5,0) -- (1.5,0); 
\draw[line width=1.2] (2.5,0) -- (5.5,0); 
\draw[line width=1.2] (6.5,0) -- (9,0); 

\node at (-9.5, 0.5){$a_0=a$}; 
\node at  (-5,1){$a_1$}; 
\node at (-2,1){$a_2$}; 
\node at (2,1){$\cdots$}; 
\node at (6,1){$a_r$}; 
\node at (10.5, 0.6){$a_{r+1}=b$}; 

\node at (-5, -1){$-$}; 
\node at (-2, -1){$+$}; 
\node at (2, -1){$-$}; 
\node at (6, -1){$-$}; 
     
   \etp  
   \caption{Transport with avoidances.}
   \label{fig:avoid}
     \end{figure}

    \paragraph{Rectilinear transport with avoidances.}  Let us now allow $[a,b]$ to contain
    other elements of $A$, say $[a,b]\cap A = \{a_0=a, a_1, \cdots, a_r, a_{r+1}=b\}$, numbered
    in the direction from $a$ to $b$ , with  $r\geq 0$.   Let $\beps = (\eps_1, \cdots, \eps_r)$, $\eps_i\in\{+, -\}$
    be a sequence formed by $+$ and $-$ signs.  We define the {\em rectilinear transport with
    avoidances} given by $\beps$ to be the map
    \be \label{eq:mabeps}
    m_{ab}^\eps = m_{ab}^{\eps, \Fc}: \Phi_{a, [a,b]} (\Fc) \lra \Phi_{b, [a,b]}(\Fc), \quad m_{ab}^\beps := m_{ab}(\gamma_\beps), 
    \ee   
    where $\gamma_\beps$ is the perturbation of the path  $[a,b]$
     obtaining by avoiding $a_i$ on the left, if $\eps_i=-$, and on the right,
     if $\eps_i=+$, see Fig. \ref{fig:avoid}. 
     
         We define $m_{ab}^+ := m_{ab}^{+, \cdots,  +}$, resp. $m_{ab}^- := m_{ab}^{-, \cdots, -}$
          to be the transport 
     with all avoidances on the right, resp. on the left (both understood as $m_{ab}[a,b]$ for $r=0$). 
     Let us note the following consequences of the Picard-Lefschetz identities   
      (Proposition \ref{prop:PL-form}). Here and later in the paper
       the identification of the stalks of the local $\Phi$-systems at the intermediate
     points is done by clockwise rotation as in  the Picard-Lefschetz 
     formula, see \eqref{eq:3-isotop}. 
     
     \begin{prop}\label{prop:PL-conseq}
     (a) We have
     \[
     m_{ab}^- \,=\,\sum_{s=0}^r  \, \sum_{1\leq i_1 <\cdots < i_s\leq r} 
     m_{a_{i_s}, b}^+ m_{a_{i_{s-1}}, a_{s_r}}^+ \cdots m_{a, a_{i_1}}^+. 
     \]
     (b) Equivalently, we have
     \[
     m_{ab}^- \,=\, m_{ab}^+ \, + \sum_{i=1}^{r} m_{a_i, b}^+ m_{a, a_i}^-.
     \]
     (c) For the composition of rectilinear transports we have the identity
     \[
     m_{a_r, b} m_{a_{r-1}, a_r} \cdots  m_{a, a_1} = 
     \sum_{\eps\in\{+,-\}^r} (-1)^{|+(\eps)|} m_{ab}^\eps,
     \]
     where $|+(\beps)|$ is the number of $+$  signs in $\beps$.
     \end{prop}
     
     \noindent{\sl Proof:} All three statements follow easily by iterated application of  
       Proposition \ref{prop:PL-form}. Part (b) is reduced to (a) by expanding each $m_{a, a_i}^-$
       according to (a). 
     Part (c) is an instance of \cite[Cor.1.1.19]{KSS} for the case of the composition
     of  paths which are rectilinear. \qed

    \paragraph{The alien derivative transport. }
    Assume now that $\ch (\k)=0$. 
   Adapting the approach of \'Ecalle, we give the following
     
     \begin{defi}\label{def:alien-tr}
     The {\em alien derivative transport} from $a$ to $b$ for $\Fc\in\Perv(\CC)$ is the map
     \[
     m_{ab}^\Delta  = m_{ab}^{\Delta, \Fc}\,= \sum_{s=0}^r {(-1)^{s+1}\over s+1} \sum_{1\leq i_1 < \cdots < i_s  \leq r}
     m_{a_{i_s}, b}^+ m_{a_{i_{s-1}}, a_{s_r}}^+ \cdots m_{a, a_{i_1}}^+: 
      \Phi_{a, [a,b]} (\Fc) \lra \Phi_{b, [a,b]}(\Fc). 
     \]
    \end{defi}

 \noindent Here the superscript $\Delta$ is just a symbol
 chosen to invoke the standard notation for alien derivatives
    \cite{ecalle, sauzin}.  
 Further, the formulas of \'Ecalle extend to the  context of perverse sheaves  in the form:      
     
     \begin{prop}\label{prop:ecalle+}
     We have
     \[
     m_{ab}^\Delta \,=\sum_{\beps\in\{+,-\}^r }
     {
    ( |+(\beps)|!) \cdot  (|-(\beps)|!)  \over (r+1)!
     }
     m_{ab}^\beps. 
     \]
     where $|+(\beps)|$ and $|-(\beps)|$ are the numbers of $+$ and $-$ signs in $\beps$
     \end{prop}

  \begin{ex}
  For $r=0$: 
  \raisebox{-0.6ex}[2ex][1ex]{
   \btp
    \draw (-1,0.5) -- (1,0.5);  
    \node at (-1, 0.5){\small$\bullet$}; 
     \node at (1, 0.5){\small$\bullet$}; 
     \node at (-1.3,0.6){$a$}; 
      \node at (1.3,0.65){$b$}; 
   \etp 
   }
   we have $m_{ab}^\Delta = m_{ab}[a,b]$. 
  
   \vskip .2cm 
   
  \noindent  For $r=1$:  \raisebox{-0.6ex}[2ex][1ex]{
   \btp
    \draw (-1,0.5) -- (1,0.5);  
    \node at (-1, 0.5){\small$\bullet$}; 
     \node at (1, 0.5){\small$\bullet$}; 
     \node at (-1.3,0.6){$a$}; 
      \node at (1.3,0.65){$b$}; 
     \node at (0,0.5){\small$\bullet$}; 
      \node at (0, 0.8){$a_1$}; 
   \etp 
   }
    we have  $m_{ab}^\Delta = {1\over 2} m_{ab}^+ + {1\over 2} m_{ab}^-$.

  \vskip .2cm
  
  For $r=2$:  \raisebox{-0.6ex}[2ex][1ex]{
   \btp
    \draw (-1,0.5) -- (1,0.5);  
    \node at (-1, 0.5){\small$\bullet$}; 
     \node at (1, 0.5){\small$\bullet$}; 
     \node at (-1.3,0.6){$a$}; 
      \node at (1.3,0.65){$b$}; 
     \node at (-0.3,0.5){\small$\bullet$}; 
      \node at (-0.3, 0.8){$a_1$}; 
        \node at (0.3,0.5){\small$\bullet$}; 
      \node at (0.3, 0.8){$a_2$}; 
   \etp 
   }
   we have
   \[
   m_{ab}^\Delta \, = \, {1\over 3} m_{ab}^{++} + {1\over 6} m_{ab}^{+-} + {1\over 6} m_{ab}^{-+} 
   + {1\over 3} m_{ab}^{--}. 
   \]
  \end{ex}
  
  \begin{rems}
  (a) Proposition \ref  {prop:ecalle+} represents $m_{ab}^\Delta$ as a linear combination of the
  $m_{ab}^\beps$ with positive coefficients summing to $1$.
  
  \vskip .2cm
  
  (b) It also shows that $m_{ab}^\Delta$ is stable under introducing dummy singularities
  and depends only on $[a,b]$ and not on $A$. 
  That is, if $\Fc$ does not really have a singularity at some $a_i$, i.e., $\Fc\in \Perv(\CC, A \, \- \, \{a_i\})
  \subset \Perv(\CC, A)$,
  then calculating $m_{ab}^\Delta$ while taking $a_i$ into account and  while not taking it into account
  gives the same answer. 
  \end{rems}
  
  \noindent {\sl Proof of Proposition \ref {prop:ecalle+}:}  This is a formal consequence
  of the Picard-Lefschetz identities (Proposition \ref{prop:PL-form}). To organize the
  calculations, let us extend the notation $m_{ab}^\eps$ to the case when
  $a=a_0, a_1,\cdots, a_r, a_{r+1}=b$ lie, in this order, on a possibly curvilinear simple
  path $\gamma$ from $a$ to $b$ which contains no other elements of $A$.
   That is, we define the path $\gamma_\eps$  as
  in \eqref{eq:mabeps}  but as the perturbation of $\gamma$, not $[a,b]$, according to $\eps$. 
  To indicate the dependence on $\gamma$, we write $m_{ab}^\eps(\gamma)$. 
  A curvilinear version of Proposition \ref{prop:PL-conseq}(c), i.e., 
 \cite[Cor.1.1.19]{KSS} gives:
  
  \begin{prop}\label{prop:vassiliev} 
  Let $\gamma_i$ denote the part of $\gamma$ between $a_i$ and $a_{i+1}$. Then
  \[
  m_{a_{r+1}, a_r}(\gamma_r) m_{a_r, a_{r-1}}(\gamma_{r-1}) \cdots
  m_{a_0, a_1}(\gamma_0) \,=\sum_{\delta\in\{+,-\}^r} (-1)^{|+(\delta)|} m_{ab}^\delta(\gamma).  \qed
  \]
  \end{prop}

  Now, to prove Proposition \ref {prop:ecalle+}, we apply Proposition \ref {prop:vassiliev} 
  to each composition in the RHS of Definition \ref {def:alien-tr}. That is, we take
  as $\gamma$ the path which goes from $a_0=a$ to $a_{i_1}$ avoiding the intermediate
  $a_i$ on the left, then from $a_{i_1}$ to $a_{i_2}$ with similar avoidances and
  continues like this, ending in the curved segment from $a_{i_s}$ to $a_{i_{s+1}}=b$
  with similar avoidances. The intermediate points are $a_{i_1}, \cdots, a_{i_s}$. 
  Then the  LHS of the formula of  Proposition \ref {prop:vassiliev}  for such $\gamma$
  and such choice of the intermetiate points is the composition in 
  Proposition \ref {prop:ecalle+}. So we get 
  \[
  \begin{gathered}
   m_{a_{i_s}, a_{i_{s+1}}}^+ m_{a_{i_{s-1}}, a_{s_r}}^+ \cdots m_{a_{i_0}, a_{i_1}}^+ \, 
   = \sum_{\delta\in\{+, -\}^s} (-1)^{|+(\delta)|} 
   m_{ab}^{(+^{{i_1}-1}, \delta_1, +^{i_2-i_1-1}, \delta_2, \cdots, \delta_s, +^{r-i_s} )},   \end{gathered} 
  \]
  where $+^m$ stands for the sequence of $m$ plus signs. 
  
  Now we need to find the coefficient at each $m_{ab}^\eps$, $\eps\in\{+,-\}^r$
  after we sum these expansions over all $s$ and all $1\leq i_1<\cdots i_s\leq r$
  with coefficients $(-1)^{s+1}/(s+1)$. For this, let us encode $\eps$ by the subset
  $I= +(\eps) = \{ i|\eps_i=+\} \subset \{1,\cdots, r\}$.  
  The coefficient is then
  \[
  \sum_{J\supset I} { (-1)^{|J\- I |}\over |J|+1 }  = \sum_{k=0}^{r-|I|} {(-1)^k\over |I|+k+1}
  { r-|I| \choose k}. 
  \]
  So Proposition \ref {prop:ecalle+} reduces to the following. 
  
  \begin{lem}
  For any integer $a,m>0$ we have
  \[
  \sum_{k=0}^m { (-1)^k\over a+k+1} {m\choose k} \,=\, {m! a! \over (m+a+1)! }. 
  \]
  \end{lem}
  
  \noindent {\sl Proof of Lemma:} For any function $f=f(a)$ of an integer variable $a$
  let $\Delta f$ be its difference derivative: $(\Delta f)(a) = f(a) - f(a+1)$. 
  The $m$th iteration of $\Delta$ has the form
  \[
  (\Delta^m f)(a) = \sum_{k=0}^m (-1)^k {m\choose k} f(a+k). 
  \]
  Let $f_m(a) = m! a!/(m+a+1)!$, $m\geq 0$. The lemma means that $f_m=\Delta^m f_0$. 
  To see this, it is enough to show that $\Delta f_m = f_{m+1}$, which is straightforward:
  \[
  \begin{gathered}
  (\Delta f_m)(a) \,=\, {m! a! \over (m+a+1)! } - {m!  (a+1)! \over (m+a+2)!} 
  \, = \, m! {a! (m+a+2)! - (a+1)! \over (m+a+2)!} \,=
  \\
  =\, m! {a! (m+1) \over (m+a+2)^!} \, = \, {a! (m+1)! \over (m+a+2)^!} \, =\, f_{m+1}(a). 
  \qed
  \end{gathered} 
  \]
  
 \paragraph{Description of $\ol\Perv(\CC,A)$ via alien transports.} We now generalize Proposition
 \ref{prop:Locperv-gen} to the case when $A = \{a_1,\cdots, a_N\} \subset \CC$ is an arbitrary finite subset. 
 We keep the notation of \S\ref {par:rect-trans} above and assume only that no $[a_i, a_j]$ is horizontal, 
 i.e,  all $\zeta_{ij}\notin\RR$, $i\neq j$. 
 
 Define the  functor  
 \[
 \Xi_\Delta: \ol \Perv(\CC,A) \lra \Mc_N, \quad \Fc\mapsto (\Phi_i, m_{ij}), 
\] 
where, as before, 
\[
\Phi_i = \bPhi_{a_i, a_i+\RR}(\Fc) = \bPhi_i(\Fc)_1, \quad  m_{ii} = \Id-T_i(\Fc)
\]
and for $i\neq j$  the map  $m_{ij}$ is the composition
\[
\bPhi_i(\Fc)_1\buildrel T_1^{\zeta_{ij}}\over  \lra \bPhi_i(\Fc)_{\zeta_{ij}} =\bPhi_{a_i,[a_i, a_j]}(\Fc)
 \buildrel m_{a_i, a_j}^{\Delta, \Fc} \over\lra
\bPhi_{a_j, [a_j, a_i]}(\Fc) = \bPhi_j(\Fc)_{\zeta_{ji}} \buildrel T_{\zeta_{ji}}^1 \over\lra  \bPhi_j(\Fc)_1. 
\]
That is, we replace the rectilinear transform in
\eqref{eq:m_{ij}-rect-GMV} (which may no longer make sense because of the presence of  intermediate points)
by  $m_{a_i, a_j}^{\Delta, \Fc}$, the  alien transform from $a_i$ to $a_j$ for $\Fc$.

\begin{prop}\label{prop:xi-delta}
The functor $\Xi_\Delta$ is an equivalence of categories. 
\end{prop}

\noindent{\sl Proof:} Let us make a small deformation of the set $A$, replacing it with $A'=\{a'_1, \cdots, a'_N\}$
with $|a_i-a'_i|\ll 1$ such that $A'$ is now in linearly general position including $\oo$. 
Since perverse sheaves are topological objects,  the continuous deformation $a_i(t) = (1-t)a_i + t a'_i$, $t\in[0,1]$
of the sets of singularities
gives rise to an equivalence $u: \ol\Perv(\CC,A) \to \ol\Perv(\CC, A')$ (``isomonodromic deformation of perverse sheaves'',
see, e.g.,  \cite{KSS}). For any $\Fc\in\ol\Perv(\CC,A)$ we denote $\Fc'$ the corresponding object of
$\ol\Perv(\CC, A')$, so that $\bPhi_i(\Fc)$ is identified with $\bPhi_{a'_i}(\Fc)$ as a local system on $S^1$. 

 Applying Proposition \ref{prop:Locperv-gen} to $A'$, we see that any $\Fc\in\ol\Perv(\CC,A)$
is uniquely determined by  the data of:

\begin{itemize}

\item[(1)] The monodromies of $\bPhi_{a'_i}(\Fc')$ which are identified with the monodromies of $\bPhi_i(\Fc)$. 

\item[(2)] The rectilinear transports for $\Fc'$.
\end{itemize} 

\noindent Now, each rectilinear transport $m_{ij}$ for $\Fc'$ corresponds, under the  equivalence $u$, to the
 rectilinear transport with
avoidances $m_{ij}^{\eps(i,j)}$ for some $\eps(i,j)$ describing to which side of $[a'_i, a'_j]$ the (formerly)  intermediate
points $a'_k$ now lie. Let $r_{ij} = |A\cap (a_i, a_j)|$ be the number of the intermediate points on $[a_i, a_j]$,
so  $\eps(i,j)$ is a sequence of length $r_{ij}$. 

The alien derivative transport $m_{a_i, a_j}^\Delta$ is a linear combination of all $2^{r_{ij}}$ transports $m_{a_i, a_j}^\eps$
with strictly positive (in particular, nonzero) coefficients. Knowing any one $m_{a_i, a_j}^\eps$, any other $m_{a_i, a_j}^{\eps'}$
is expressed, in virtue of   the Picard-Lefschetz formulas, by adding or subtracting compositions of transports with avoidances for
smaller subintervals of $[a_i, a_j]$. Hence the data of  all $\{m_{a_i, a_j}^\eps\}$ for all distinct  $1\leq i,j \leq N$
and all $\eps\in \{+, -\}^{r_{ij}}$,  is uniquely
recovered (by triangular-type formulas) from the data of $\{m_{a_i, a_j}^{\eps(i,j)}\}$, $1\leq i,j\leq N$,
where we choose one representative $\eps(i,j)$ for each ordered pair $(i,j)$. Therefore the data of
such  $\{m_{a_i, a_j}^{\eps(i,j)}\}$ are in bijection with the data of  $\{m_{a_i, a_j}^\Delta\}$, and the proposition 
is proved. \qed
  

  \subsection{Alien derivatives and Stokes automorphisms 
   for perverse sheaves}\label{subsec:alien-der}
  
  \paragraph{Multiplicative properties of one-sided avoidances.} Again, we start by
  taking $\k$ to be an arbitrary field. 
  For $a,b\in\CC$ and a perverse sheaf $\Fc\in\Perv(\CC)$ we use notation
  $m_{ab}^\pm = m_{ab}^{\pm, \Fc}$ to mean either  $m_{ab}^{+, \Fc} = m_{ab}^{+,\cdots, +, \Fc}$
  or $m_{ab}^{-, \Fc} = m_{ab}^{-,\cdots, -, \Fc}$, this meaning  to be used
  consistently in any formula. 
  
  \vskip .2cm
  
    Let  $\Fc \in\Perv^0(\CC, A')$, $\Gc\in\Perv^0(\CC,A'')$ and  $A=A'+A''$. 
  By Theorem \ref {thm:Thom},  for any direction $\theta\in S^1$ and any $a,b\in A$
  we have 
    \be\label{eq:thom-thom}
  \bPhi_a(\Fc * \Gc)_\theta \= \bigoplus_{a'\in A', a''\in A'' \atop
  a'+a''=a} \bPhi_{a'}(\Fc)_\theta \otimes\bPhi_{a''}(\Gc)_\theta, \quad 
    \bPhi_b(\Fc * \Gc)_{-\theta} \= \bigoplus_{b'\in A', b''\in A''\atop b'+b''=b}
     \bPhi_{b'}(\Fc)_{-\theta} \otimes\bPhi_{b''}(\Gc)_{-\theta}.
  \ee
  Let $a = a'+a'', b=b'+b''\in A$ be distinct, with $a',b'\in A'$ and $a'', b''\in A''$. 
  Note that it is possible that $a'=a''$ or $b'=b''$ (but not both). 
  Take  $\theta = \zeta_{ab}$ to be the direction from $a$ to $b$. Let us view the
  rectilinear transport with one-sided avoidances as a linear map
  \[
  m_{ab}^{\pm, \Fc * \Gc} : \bPhi_a(\Fc*\Gc)_{\theta} \lra \bPhi_b (\Fc* \Gc)_{-\theta}.
  \]
   With respect
  to the decompositions \eqref{eq:thom-thom}, we then have the matrix element
  \[
  \bigl(m_{ab}^{\pm, \Fc * \Gc}\bigr)_{a',a''}^{b', b''}: \bPhi_{a'}(\Fc)_\theta \otimes\bPhi_{a''}(\Gc)_\theta 
  \lra
  \bPhi_{b'}(\Fc)_{-\theta} \otimes\bPhi_{b''}(\Gc)_{-\theta}.
  \]

  \begin{thm}\label{thm:mult-mat-el}
  (a) Unless the intervals $[a',b'], [a'', b'']$ and $[a,b]$ are parallel with the same direction,
  $ \bigl(m_{ab}^{\pm, \Fc * \Gc}\bigr)_{a',a''}^{b', b''}=0$. Here 
   a degenerate interval $[a',a']$ or $[b',b']$  is considered parallel  (with the same 
   direction)
   to any other interval.
   
   \vskip .2cm
   
   (b) If  the intervals $[a',b'], [a'', b'']$ and $[a,b]$ are parallel with the same direction, then 
   \[
    \bigl(m_{ab}^{\pm, \Fc * \Gc}\bigr)_{a',a''}^{b', b''}=m_{a',b'}^{\pm, \Fc} \otimes m_{a'', b''}^{\pm, \Gc}.
    \]
    Here we understand $m^{\pm, \Fc}_{a',a'}$ or $m^{\pm, \Gc}_{a'', a''}$ as 
    the identity map. 
  
  \end{thm}
  
  \noindent Let us express the above 
  condition of three intervals being parallel with the same direction
  by $  [a', b'] \parallel  [a'', b''] \parallel [a,b]$. Then we can reformulate
  Theorem \ref {thm:mult-mat-el} as follows:
  
 \begin{refo}\label{ref:m-pm-mult}
 In the above notation, we have 
    \[
  m_{ab}^{\pm, \Fc  * \Gc} \,= \sum_{a'+a''=a, \, b'+ b''= b
  \atop
  [a', b'] \parallel  [a'', b''] \parallel [a,b]] } 
  m_{a',b'}^{\pm, \Fc} \otimes m_{a'', b''}^{\pm, \Gc}. \qed
  \]
  
  \end{refo}

 \noindent Theorem \ref {thm:mult-mat-el} and Reformulation \ref {ref:m-pm-mult}  as well as the proof
  below
  are inspired by Theorem 6.83 of \cite{sauzin} and its purely
  analytic proof.

  \paragraph{ Proof of Theorem \ref {thm:mult-mat-el}.} 
  
\begin{figure}[h]
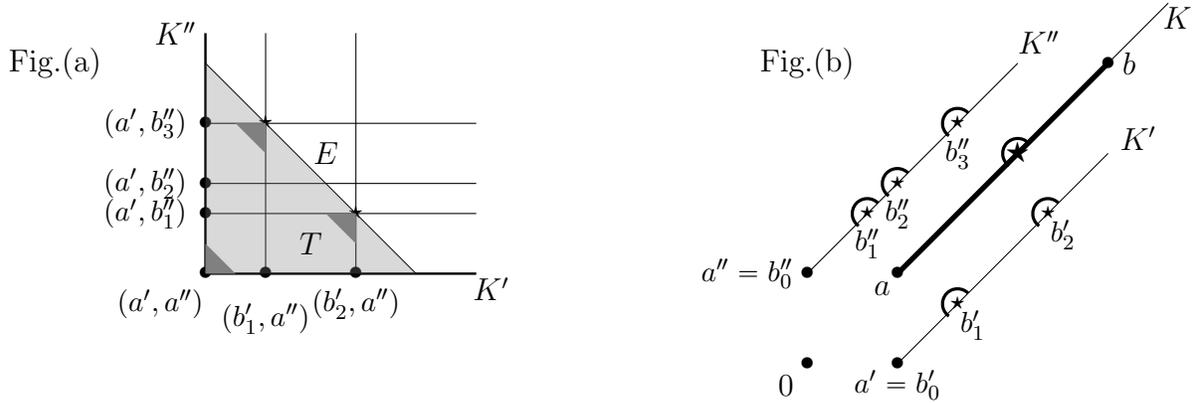

     \centering
\btp[scale=.4, baseline=(current  bounding  box.center)]

\node at (0,10){Fig.(b)}; 

\node at (0,0) {$\bullet$}; 
\node (a')  at (0,3) {$\bullet$}; 
\node  (a'') at (3,0) {$\bullet$}; 
\node  (a)  at (3,3) {$\bullet$}; 
\node (b)  at (10, 10){$\bullet$}; 
\node (b'_1) at (5,2) {$\star$};
\node (b'_2) at (8,5) {$\star$}; 
\node (b"_1) at  (2,5) {$\star$}; 
\node (b"_2) at (3,6) {$\star$}; 
\node (b"_3) at (5, 8) {$\star$}; 
\node (b'_1+b''_1) at (7,7) {$\bigstar$}; 

\draw (3,0) -- (10,7); 
\draw (0,3) -- (7,10); 
\draw (3,3) -- (12,12); 
\draw [line width = 2] (3,3) -- (10,10); 

\draw[line width=1.2] (5.35, 2.35) arc [start angle=45, end angle=235, x radius=0.5cm, y radius=0.5cm];
\draw[line width=1.2] (8.35, 5.35) arc [start angle=45, end angle=235, x radius=0.5cm, y radius=0.5cm];
\draw[line width=1.2] (2.35, 5.35) arc [start angle=45, end angle=235, x radius=0.5cm, y radius=0.5cm];
\draw[line width=1.2] (3.35, 6.35) arc [start angle=45, end angle=235, x radius=0.5cm, y radius=0.5cm];
\draw[line width=1.2] (5.35, 8.35) arc [start angle=45, end angle=235, x radius=0.5cm, y radius=0.5cm];
\draw[line width=1.2] (7.35, 7.35) arc [start angle=45, end angle=235, x radius=0.5cm, y radius=0.5cm];

\node at (-0.7, -0.7) {$0$}; 
\node at (3, -0.7) {\small $a'=b'_0$}; 
\node at (5.5, 1.3) {\small $b'_1$}; 
\node at (8.5, 4.3) {\small $b'_2$}; 

\node at (11,7.5) {$K'$}; 
\node at (12.3, 11.5){$K$}; 

\node at (2.5, 2.5) {$a$}; 
\node at (10.7, 10){$b$}; 
\node at (-2,3){\small $a''=b''_0$}; 
\node at  (2,4) {\small $b''_1$}; 
\node at  (3,5) {\small $b''_2$}; 
\node at (5,7) {\small $b''_3$}; 

\node at (7.7,10.7) {$K''$}; 


\node at (-25, 10) {Fig.(a)}; 

\node at (-20, 3){$\bullet$}; 
\draw [line width = 0.9] (-20,3) -- (-11, 3); 
\draw [line width = 0.9] (-20,3) -- (-20, 11); 
\node at (-18,3) {$\bullet$}; 
\node at (-15,3) {$\bullet$}; 
\node at (-20, 5) {$\bullet$}; 
\node at (-20,6){$\bullet$}; 
\node at (-20,8) {$\bullet$}; 

\draw (-20, 10) -- (-13,3); 

\draw (-18,3) -- (-18,11); 
\draw (-15,3) -- (-15, 11); 
\draw (-20, 5) -- (-11,5); 
\draw (-20,6) -- (-11,6); 
\draw (-20,8) -- (-11,8); 

\node at (-18,8) {$\star$}; 
\fill [gray] (-18,8) -- (-18,7)  -- (-19,8) -- (-18,8); 

\node at (-15,5) {$\star$}; 
\fill [gray] (-15,5) -- (-15,4) -- (-16,5) -- (-15,5);

\node at (-10.5, 2.5) {$K'$}; 
\node at (-21,11){$K''$}; 
\node at (-21.5,2){\small$(a',a'')$}; 

\node at (-18,1.5){\small $(b'_1, a'')$}; 

\node at (-15,2) {\small $(b'_2, a'')$}; 
\node at (-22,5){\small $(a', b''_1)$}; 
\node at (-22,6){\small $(a', b''_2)$}; 
\node at (-22, 8){\small $(a', b''_3)$}; 

\fill[gray, opacity=0.3] (-20,3) -- (-20,10) -- (-13,3) -- (-20,3); 

\node at (-16.5,4) {$T$}; 

\fill  [gray, opacity=1] (-20,3) -- (-20,4) -- (-19,3) -- (-20,3); 
\node at (-16, 7) {$E$}; 
\etp
\caption{The area $T\subset K'\times K''$ and  the transport with avoidances for $\Fc * \Gc$ . }
\label{fig:convo}
\end{figure}

  Let us treat the case of $m^-$ (avoidances on the left), the case of $m^+$ being
  similar. 
   Fix  $a'\in A'$, $a''\in A''$.  Let $a=a'+a''\in A$ and
   \be\label{eq:K'K''}
   \begin{gathered}
 K'= a'+\RR_+ \theta = \{ a' + t' \, \theta |\,\, t' \geq 0\}, \quad  K'' =  a''+\RR_+  \theta
  = \{a''+ t'' \, \theta | \,\, t'' \geq 0\},
  \\
 K= a+ R_+\theta = \{a+t\,\theta| \,\, t\geq 0\}
 \end{gathered}
   \ee
  be the straight half-lines issuing from $a'$,  $a''$ and $a$ in the direction $\theta$,
  see  Fig. \ref{fig:convo}(b). 
  We use $t'$,  $t''$  and $t$ as coordinates on these half-lines. 
  
      By
  Example \ref {ex:Phi-Psi},
 \be\label{eq:Phi-Phi}
 \begin{gathered}
 \bPhi_{a'}(\Fc)_\theta \,=\, \ul H^0_{K'}(\Fc)_{a'}, 
 \quad \bPhi_{a''}(\Gc)_\theta \,=\, \ul  H^0_{K''}(\Fc)_{a''} \quad \text{and therefore}
 \\
  \bPhi_{a'}(\Fc)_\theta \otimes \bPhi_{a''}(\Gc)_\theta \,=\, \Hc_{(a', a'')}, \quad
  \text{where} \quad \Hc := 
  \ul H^0_{K'\times K''}
  (\Fc \boxtimes\Gc).
 \end{gathered}
 \ee
 The stalk $\Hc_{(a', a'')}$ can be seen as the $0$th cohomology of $\Fc\boxtimes\Gc$
 with support in the dark shaded area near the left of $(a', a'')$ in 
 Fig. \ref{fig:convo}(a). 
    Denote 
    \[
    \wt A{'} = A'\cap K' \, =\, \bigl\{a'=b'_0, b'_1, b'_2, \cdots\bigr\}, \quad
    \wt A{''}= A''\cap K'' \, = \, \bigr\{ a''=b''_0, b''_1, b''_2, \cdots \bigr\}
    \] 
    in order given by the direction of $K', K''$, see Fig. \ref{fig:convo}(b).

The sheaf  $\Hc$  on $K'\times K''$ is constructible with respect to
the stratification  cut out by  $\wt A{'}\times K''$,
    $K'\times \wt A{''}$ and their intersection $\wt A{'}  \times \wt A{''}$,
    see  Fig. \ref{fig:convo}(a). This is because $\Fc\boxtimes\Gc$ is
    constructible w.r.t. a similar stratification of $\CC\times\CC$. 
    In particular, $\Hc$ is locally constant on the
   interior of $K'\times K''$ near $(a', a'')$.
   
   \vskip .2cm
   
   Let $\wt K{'}$, $\wt K{''}$ be small perturbations of $K', K''$ obtained by avoiding 
   the $b'_i, b''_i, i>0$, on the left. Then $\wt K{'} \times \wt K{''}$
   coincides with $K'\times K''$ near $(a', a'')$. By construction, 
   \[
   \wt\Hc \, = \, \ul H^0_{\wt K{'}\times\wt K{''}} (\Fc\boxtimes\Gc)
   \]
   is locally constant on the entire interior of $\wt K{'} \times \wt K{''}$ and coincides
   with $\Hc$ near $(a', a'')$.  Let also $\wt K$ be a  similar perturbation of 
   $[a,b]\subset K$ avoiding
   all the elements of $A$ other than $a$ and $b$ on the left. 
   
   \vskip .2cm
   
    Now look
  at the composite map
  \be\label{eq:mat-elem1}
  \bPhi_{a'}(\Fc)_\theta  \otimes \bPhi_{a''}(\Gc)_\theta
  \buildrel \eps_{a', a''} \over  \lra \bPhi_a(\Fc * \Gc)_\theta 
  \buildrel m_{ab}^{-, \Fc * \Gc} \over\lra \bPhi_b(\Fc * \Gc)_{-\theta} \,=\, \bigoplus_{b'+b''=b} \bPhi_{b'} (\Fc)_{-\theta}
  \otimes \bPhi_{b''}(\Gc)_{-\theta}.
  \ee
   Here $\eps_{a', a''}$ is an embedding given by the Thom-Sebastiani theorem 
    \ref{thm:Thom} which is also indicated in the equality on the right. 
  The transport $m_{ab}^{-, \Fc*\Gc}$ is defined using $\wt K$. 
  As explained in \S  \ref{subsec:perRiem} \ref {par:transport} (applied to $\alpha = \wt K$)
 it is composed of three maps:
    
    \begin{itemize}
    \item[(1)] The generalization map $u_{a,\wt K} = u_{a,\wt K}^{\Fc * \Gc}$ (in the notation of \eqref{eq:u-a,alpha})
    from  $\bPhi_a(\Fc*\Gc)_\theta = \ul H^0_{\wt K}(\Fc* \Gc)_a$ to the stalk
     at a nearby point $c\in \wt K$ which is
    $\ul H^0_{\wt K}(\Fc * \Gc)_c =  (\Fc * \Gc)[-1]_c$, the same as the stalk at $c$ of the
    local system $\Fc * \Gc [-1]$. 
    
    \item[(2)] The parallel transport of the result of (1) along $\wt K $ in the local system 
    $\Fc * \Gc [-1]$ until we almost reach  $b$. 
    
    \item[(3)] After approaching close to $b$ using (2), 
    applying the variation map $v_{b,\wt K}= v_{b,\wt K}^{\Fc * \Gc}$  (in the notation of
     \eqref{eq:u-a,alpha})
     at $b$ which is the dual  of the generalization 
     map $u^{(\Fc * \Gc)^\vee}_{b, \wt K}$ for the 
     Verdier dual perverse sheaf $(\Fc * \Gc)^\vee \= \Fc^\vee * \Gc^\vee$. 

  \end{itemize} 
  
 \noindent  Let $\tau\in K'\times K''$ be an interior point close to 
 the point $(a', a'')$. It also lies
  in $\wt K{'}\times\wt K{''}$ and
  \[
  \gamma_{a',a''} : \Hc_{(a', a'')} = \wt\Hc_{(a', a'')} \lra\Hc_\tau = \wt\Hc_\tau
  \]
   be the generalization map of the constructible sheaves $\Hc, \wt\Hc$
   which coincide in the area containing $(a', a'')$ and $\tau$.  
   
   \vskip .2cm
   
   Let $\phi\in  \bPhi_{a'}(\Fc)_\theta  \otimes \bPhi_{a''}(\Gc)_\theta
=\Hc_{(a',a'')}$. 
   As follows from the construction of the identification in the Thom-Sebastiani
 theorem (proof of  Theorem \ref{thm:Thom}), 
   the composition  $u_{a, \wt K} \eps_{a',a''}$  can be seen as the composition
   of $\gamma_{a', a''}$ followed by the ``forgetting of support'' morphism 
   \[
   R(+)_* \,\, \ul H^0_{K'\times K''}(\Fc\boxtimes \Gc) \lra R(+)_* (\Fc\boxtimes \Gc)  = \Fc * \Gc
   \]
  evaluated in the stalks over $c$. Therefore we can replace parallel transport
  of   $u_{a, \wt K} \eps_{a',a''}(\phi)$ along $\wt K$ by  parallel transport
  of $\gamma_{a', a''}(\phi)$ in the local system given by  $\wt \Hc$    
  on the interior of $\wt K{'}\times\wt K{''}$. As this interior is contractible,
  we have a well defined section $\wt\phi$ of $\wt \Hc$ on it, extending $\gamma_{a', a''}(\phi)$.
  
  \vskip .2cm
  
  The map $+: \CC\times\CC\to\CC$ restricts to $+_K: K'\times K''\to K$
  which in coordinates $t', t''$ from \eqref{eq:K'K''} has the form $t=t'+t''$. 
  Let $T = +_K^{-1}([a,b))\subset K'\times K'' $ be the preimage of the half-open
  interval $[a,b)$,  depicted as the large shaded area on  Fig. \ref{fig:convo}(a),
    and  $\wt T\subset \wt K{'}\times \wt K{''}$ be 
  corresponding perturbation of $T$. 
  On the edge  $E=\{t'+t''=|b-a|\}$ of $T$ we have the points $(b'_i, b''_j)$ with $b'_i + b''_j=b$,
  i.e. precisely the points $(b', b'')$ such that $[a',b']$,  $[a', b'']$ and $[a,b]$
  are parallel in the same direction, as
  in the statement of Theorem  \ref {thm:mult-mat-el}.
  
  \vskip .2cm
  
Since the variation map $v_{b,\wt K}^{\Fc * \Gc}$
  in (3) above is dual to $u_{b, \wt K}^{\Fc^\vee * \Gc^\vee}$ which has been
  just described above, we see that only $(b'_i, b''_j)$ on $E$
  will receive a component of $m_{ab}^{-, \Fc}(\eps_{a,a'}(\phi))$:
  the other $(b_i', b''_j)$ will map far from $b$. 
  This is precisely part (a) of the theorem. 
  
  \vskip .2cm
  
  Let us now  prove  (b). Let $b'=b'_i, b''= b''_j$ be such that $b'+b''=b$
  and $\psi\in \bPhi_{b'}(\Fc^\vee)_{-\theta}  \otimes\bPhi_{b''}(\Gc^\vee)_{-\theta}$. 
 As before, we use the duality between $v_{b,\wt K}^{\Fc * \Gc}$
  and $u_{b, \wt K}^{\Fc^\vee * \Gc^\vee}$ and interpret  the latter
  in terms of  the generalization map $\gamma_{(b', b'')}: \Hc_{(b', b'')}\to \Hc_\sigma$
  where $\sigma$ is a point of $T$ in the area  near $(b', b'')$, depicted as the darker 
  shaded area  near the edge $E$ on Fig. \ref{fig:convo}(b). So we have the equality of the
  pairings
  \[
  \bigl( m_{ab}^{-, \Fc * \Gc}(\eps_{a', a''}(\phi)), \psi\bigr) =
  \bigl(\wt\phi(\sigma) , \gamma_{(b', b'')}(\psi)\bigr). 
  \]
  where on the right the sections of the two dual local systems are evaluated
  at a nearby points so the pairing is well defined. 
  But since $\wt\phi$ is a section of the local system
  $\Fc\boxtimes\Gc[-2]$,  pairing on the right is precisely 
  $\bigl( (m_{a',b'}^{-,\Fc} \otimes m_{a'', b''}^{-, \Gc}) (\phi), \psi\bigr)$. 
 This proves the theorem. 
 
\begin{figure}[h]
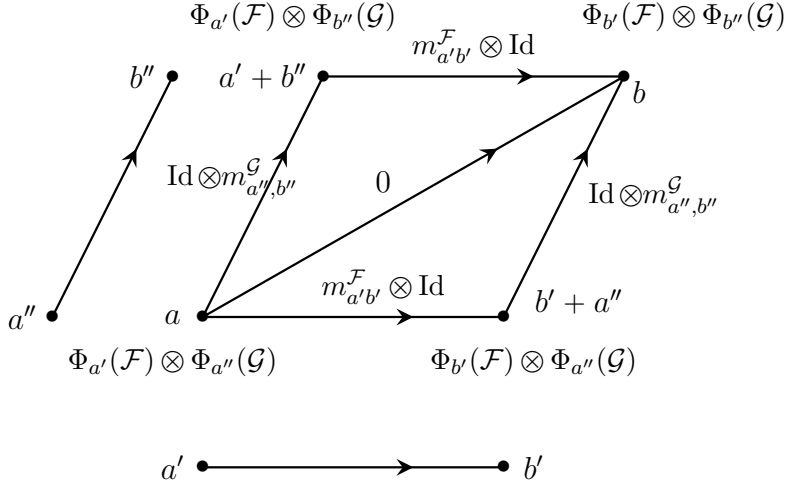

     \centering
\btp[scale=.4, baseline=(current  bounding  box.center)]

\node (a') at (5,0){$\bullet$}; 
\node at (4,0){$a'$};

\node (a'') at (0,5){$\bullet$}; 
\node at ((-1,5) {$a''$}; 
\node (a)  at (5,5){$\bullet$}; 
\node at (4,5){$a$}; 
\node at (4 ,3.5) {\small $\Phi_{a'}(\Fc) \otimes\Phi_{a''}(\Gc)$}; 

\node (b'') at (4,13){$\bullet$}; 
\node at (3,13) {$b''$}; 
\node (b') at (15,0){$\bullet$}; 
\node  at (16,0) {$b'$}; 
\node (a'+b'') at (9,13){$\bullet$}; 
\node at (7,13) {$a'+b''$}; 
 \node at (6, 9.6) {\small $\Id \otimes m_{a'', b''}^\Gc$}; 

\node at (8 ,15) {\small $\Phi_{a'}(\Fc) \otimes\Phi_{b''}(\Gc)$}; 
\node at (14,  14) {\small $m_{a' b'}^\Fc \otimes\Id$}; 

\node (b'+a'') at (15,5){$\bullet$}; 
\node at (17.5, 5.5) {$b'+a''$}; 
\node at (16, 3.5) {\small $\Phi_{b'}(\Fc) \otimes\Phi_{a''}(\Gc)$}; 
 \node at (20, 9) {\small $\Id \otimes m_{a'', b''}^\Gc$}; 

\node (b) at (19,13) {$\bullet$}; 
\node at (19.5, 12.5) {$b$}; 
\node at (21, 15) {\small$\Phi_{b'}(\Fc) \otimes\Phi_{b''}(\Gc)$}; 

\draw   [ decoration={markings,mark=at position 0.7 with
{\arrow[scale=1.5,>=stealth]{>}}},postaction={decorate},
line width = .3mm] (5,0) -- (15,0) ;  

\node at (11,6) {\small $m_{a' b'}^\Fc \otimes\Id$}; 

\node at (11, 9.5) {$0$}; 

\draw   [ decoration={markings,mark=at position 0.7 with
{\arrow[scale=1.5,>=stealth]{>}}},postaction={decorate},
line width = .3mm] (0,5) -- (4,13) ;  

\draw  [ decoration={markings,mark=at position 0.7 with
{\arrow[scale=1.5,>=stealth]{>}}},postaction={decorate},
line width = .3mm] (5,5) -- (15,5) ;  

\draw  [ decoration={markings,mark=at position 0.7 with
{\arrow[scale=1.5,>=stealth]{>}}},postaction={decorate},
line width = .3mm] (5,5) -- (9,13) ;  

\draw  [ decoration={markings,mark=at position 0.7 with
{\arrow[scale=1.5,>=stealth]{>}}},postaction={decorate},
line width = .3mm] (9,13) -- (19,13) ;  

\draw   [ decoration={markings,mark=at position 0.7 with
{\arrow[scale=1.5,>=stealth]{>}}},postaction={decorate},
line width = .3mm]  (15,5) -- (19,13) ;

\draw [ decoration={markings,mark=at position 0.7 with
{\arrow[scale=1.5,>=stealth]{>}}},postaction={decorate},
line width = .3mm] (5,5) -- (19,13) ;

\etp
\caption{An elementary parallelogram: the diagonal transport is $0$. }
\label{fig:paral}
\end{figure}

 \paragraph{Example: an elementary parallelogram.} As an illustration of Theorem
  \ref {thm:mult-mat-el} consider the following particular case.
  Let $A'= \{a',b''\} $ and $A'' = \{ a'', b''\}$  each consist of two elements
   such that the intervals 
  $[a',b']$ and $[a'', b'']$
  are not parallel, so    
    \[
  A = A'+ A'' \,=\, \bigl\{ a:=a'+a'', \,\,b'+a'', \,\, a'+ b'', \,\, b'+b''=: b \bigr\} 
  \]
  is the set of vertices of  a nondegenerate parallelogram,  see Fig. \ref{fig:paral}.
 Then by Theorem \ref{thm:Thom} the vanishing cycle spaces of $\Fc*\Gc$ at these vertices
 are the tensor products, as indicated in  Fig. \ref{fig:paral}.

\noindent  As there are
 no intermediate points, 
the rectilinear transports for $\Fc * \Gc$
 between these vertices do not need avoidances: $m^+ = m^- = m$.
 In this situation, 
 Theorem   \ref {thm:mult-mat-el}  says that the transports along the
 faces of the parallelogram are tensor products of $m_{a',b'}^\Fc$
 or $m_{a'', b''}^\Gc$ with $\Id$, so these maps look (up to isomorphisms
 of the stalks of the  local systems $\bPhi$) as forming 
 a commutative square.
 
 \vskip .3cm
 
  But the diagonal transport
 $m_{ab}^{\Fc * \Gc}$ is equal to $0$.  This last statement 
 can be seen directly by noticing
 that the rays $K', K''$ coming from $a'$ and $a''$ in the direction
 $\theta=\zeta_{ab}$ will contain no other elements of $A'$ or $A''$. 
So $\wt K{'}= K'$, $\wt K{''}=K''$ and the sheaf $\Hc = \wt\Hc$ will be locally
constant everywhere inside $K'\times K''$. This means that for 
$\phi\in\Phi_{a'}(\Fc)\otimes\Phi_{a''}(\Gc)$ the image $u_{a, K}(\phi)$
of $\phi$ under the generalization map, can be continued along $K$
all the way  through $b$
and so its variation at $b$ (image under $v_{b,K}$) is zero. 

  \vfill\eject

  \paragraph{Matrix formulation. The Stokes operator. } 
   As before, let $S^1 = \{\zeta\in \CC: \, |\zeta|=1\}$ be the circle of directions. 
  Let $\omega\in\CC$ be a nonzero number. For any 
  $\Fc\in\Perv(\CC)$ we write $\bPhi_a(\Fc)_\omega$
  for the stalk of $\bPhi_a(\Fc)$ at $\omega/|\omega| \in S^1$ and put
  $\bPhi(\Fc)_\omega = \bigoplus_{a\in\CC}\bPhi_a(\Fc)_\omega$.
  
  \vskip .2cm

   Define the  operator
  $C_\omega^{\pm, \Fc}: \bPhi(\Fc)_\omega \to\bPhi(\Fc)_\omega$
   by defining its matrix
  elements 
  $(C_\omega^{\pm, \Fc})_a^b: \bPhi_a(\Fc)_\omega\to \Phi_b(\Fc)_\omega$ 
  as follows: 
  \[
  (C_\omega^{\pm, \Fc})_a^b = \begin{cases}
  \Id, & \text{ if } a=b;
  \\
  T_{-\omega}^\omega \circ m_{ab}^{\pm, \Fc}, & \text{ if } b=a+\omega;
  \\
  0, & b\text{ otherwise}.
  \end{cases}
  \]
     Here $T_{-\omega}^\omega$ is the clockwise half-monodromy of
  $\bPhi_b(\Fc)$ from the direction $-\omega/|\omega|$
  to $+\omega/|\omega|$ (same identification as used in
  composing rectilinear transports).

  Reformulation \ref {ref:m-pm-mult}  can be further reformulated  as follows.
  
  \begin{cor} \label{prop:C-omega} 
  For $\Fc, \Gc\in\Perv^0(\CC)$ we have
  \[
  C_\omega^{\pm, \Fc * \Gc} \,=\sum_{\omega'+\omega'' = \omega
  \atop
  \omega',  \omega'' \in [0, \omega] } C_{\omega'}^{\pm, \Fc} \otimes
  C_{\omega''}^{\pm, \Gc}.
  \]
Here for $\omega'=0$ or $\omega''=0$  (only one case can
occur, as $\omega\neq 0$) we understand $C_0^\pm$ 
   to be $\Id$.  \qed
  \end{cor}

  \vskip .5cm
   
 \begin{defi}\label{defi:stokes}
 Let $\zeta\in S^1$ and $\Fc\in\Perv(\CC)$. We define the   {\em Stokes operator} associated to
 $\zeta$ and $\Fc$ as
 \[
 \St_\zeta = \St_\zeta^\Fc =  \Id + \sum_{\omega\in\RR_{>0}\zeta} C^-_\omega: \bPhi(\Fc)_\zeta \lra \bPhi(\Fc)_\zeta.
 \]
 \end{defi}

 The operator $\St_\zeta$ is invertible because it is represented by a block-upper triangular matrix
 with respect to
 the order $\leq_\zeta$ on $A$.  It has $\Id$ on the diagonals since it gives identity on the
 associated graded space. 
  Proposition \ref{prop:C-omega}   can be reformulated even more concisely. 
 
 \begin {prop}\label{prop:stokes=auto}
 Let $\Fc, \Gc\in \Perv^0(\CC)$.  For any $\zeta\in S^1$. 
 \[
 \St_\zeta^{\Fc * \Gc} \,=\, \St_\zeta^\Fc \otimes\St_\zeta^\Gc.  \qed
 \]

 \end{prop}
 
 \noindent In other words, $\St_\zeta$ is an automorphism of the tensor functor $\bPhi(-)_\zeta$.

  \paragraph{Alien derivatives  via matrix elements.}
  Assume now that $\ch(\k)=0$. 
 For $\Fc\in \Perv(\CC)$ and a nonzero  
  $\omega\in\CC$  we  call the
   {\em alien derivative for $\Fc$  in the direction} $\omega$ the operator
  $\Delta_\omega =
   \Delta_\omega^\Fc: \bPhi(\Fc)_\omega\to\bPhi(\Fc)_\omega$
  whose matrix elements
  $
  (\Delta_\omega)_a ^b: \bPhi_a(\Fc)_\omega\to\bPhi_b(\Fc)_\omega
  $ are defined as follows:
  \[
  (\Delta_\omega)_a^b = \begin{cases}  T_{-\omega}^\omega \circ
  m_{ab}^\Delta, &  \text{ if } b=a+\omega;
  \\
  0,& \text{ otherwise.} 
  \end{cases}
  \]
  Thus $\Delta_\omega = 0$ for almost all $\omega$. 
  
  \paragraph{Alien derivatives as functor derivations.}
  As $\St^\Fc_\zeta$ is given by  a block-upper triangular matrix with $\Id$ on the diagonal,
  its logarithm is  a well defined operator. 
  
  \begin{thm}
  (a) We have 
  \[
  \log \St_\zeta^\Fc = \sum_{\omega \in \RR_{>0}\zeta} \Delta_\omega^\Fc.
  \]
  (b) Let $\Fc, \Gc\in\Perv^0(\CC)$.
  With respect to the identification $\bPhi(\Fc* \Gc)_\omega \=
   \Phi(\Fc)_\omega \otimes\Phi(\Gc)_\omega$
  we have the Leibniz rule
  \[
  \Delta_\omega^{\Fc*\Gc} \,=\,\Delta_\omega^\Fc \otimes \Id \, + \, \Id\otimes \Delta_\omega^\Gc. 
  \]
  
  \end{thm}
  
 \noindent  In other words, the alien derivative is a derivation of the tensor functor
   $\bPhi(-)_\omega$. 
  
  \vskip .2cm
  
  \noindent{\sl Proof:}  (a) follows by comparison of Definition \ref {def:alien-tr} of
  the $m_{ab}^\Delta$ with the logarithmic series 
  $\log(1+x) = \sum_{s=0}^\oo (-1)^s x^{s+1}/(s+1)$. 
  
  \vskip .2cm
  
  (b) Since $\St_\zeta$ is an automorphism of the  tensor functor 
  $\bPhi(-)_\zeta$, its logarithm  $\Delta_{\RR_{>0}\zeta} = \log \St_\zeta$ is a derivation by formal reasons.  Now, $\bPhi(-)_\zeta$ takes values in the tensor category of $\CC$-graded
  vector spaces (with the graded tensor product). Any endomorphism $D$ of this functor
  can be split into homogeneous components $D=\sum_{\omega\in\CC}  D_\omega$,  where
  $D_\omega$ raises the degree by $\omega$. Clearly, $D$ is a derivation if and only if
  each $D_\omega$ is a derivation.  It remains to notice that the $\Delta_\omega$,
  $\omega\in \RR_{>0}\zeta$ are precisely the homogeneous components
  of $\Delta_{\RR_{>0}\zeta}$,  in virtue of (a). \qed

\paragraph{Stokes automorphisms in terms of the Fourier transform.} Let
$\k=\CC$ and  $\Fc\in Perv(\CC, A)$. 
The directions $\zeta_{ab}$ for all distinct $a,b\in A$, see \eqref{eq:zeta-ab},
 will be called
{\em Stokes directions} for $A$.

As in \S \ref {subsec:FT},
 Fourier transform gives a local system $\FT_\gen(\Fc)$ on $\CC^*$ or,
equivalently, on  $S^1$. 
In  the proof of Proposition \ref{prop:FTgen} we constructed an identification 
of $\FT_\gen(\Fc)$
with $\bPhi(\Fc) = \bigoplus_{a\in A} \bPhi_a(\Fc)$ outsides of the Stokes directions. 
Indeed,  $\zeta\in S^1$ is non-Stokes if and only if  all the half-rays
$K_a(\zeta) = a+\zeta\RR_+$ are disjoint.

\noindent Let us now complete that construction by
describing how these identifications  glue together
at a given Stokes direction  $\zeta$. Let $\zeta^+$ and $\zeta^-$ be
nearby non-Stokes directions clockwise and anti-clockwise from $\zeta$. 
The gluing along $\zeta$ for the local system $\FT_\gen(\Fc)$).
 with respect to our prior identifications  is given by the  map  $S_\zeta$
 defined as the composition 
 \[
\bPhi(\Fc)_\zeta \buildrel T_{\zeta}^{\zeta^+}(\bPhi)\over \lra \bPhi(\Fc)_{\zeta^+} 
\buildrel {\eqref {prop:FTgen}} \over
\lra \FT_\gen(\Fc)_{\zeta^+}
\buildrel T_{\zeta^+}^{\zeta^-}(\FT) \over
\lra\FT_\gen(\Fc)_{\zeta^-} 
\buildrel  {\eqref {prop:FTgen}} \over
\lra \bPhi(\Fc)_{\zeta^-} 
\buildrel T_{\zeta^-}^\zeta(\bPhi)\over
\lra \bPhi(\Fc)_\zeta. 
\]
Here,  say,  $T_{\zeta}^{\zeta^+}(\bPhi)$ is the monodromy from $\zeta$ to $\zeta^+$
(along the shortest path) for the local system $\bPhi(\Fc)$, and
 $T_{\zeta^+}^{\zeta^-}(\FT)$ is the monodromy from $\zeta^+$ to $\zeta^-$
 for $\FT_\gen(\Fc)$. 
The following result   up to notation coincides with  \cite[Th.5.2.2]{dagnolo-sabbah}.  
 
\begin{thm}\label{thm:S-zeta-stokes}
$S_\zeta$ coincides with the Stokes automorphism $\St_\zeta$. \qed
\end{thm}

\begin{rem}
Because of the triangular nature of $\St_\zeta$, it preserves the Stokes
filtration which was described in Proposition \ref {prop:stokes-filt} for generic
(non-Stokes) directions $\theta$:
the terms that are added after crossing $\zeta$,
 have lower rate of exponential growth.
 So we obtain  a well defined filtration
on the local system $\FT_\gen(\Fc)$ on $S^1$  labelled by the sheaf of posets
$(\ul A, \leq_\theta)_{\theta\in S^1}$  on $S^1$.
This means that the sheaf of sets $\ul A$ is constant,
but the order $\leq_\theta$ varies with $\theta$, see \cite{deligne-var, deligne-docs}
and  \cite[\S2.4]{KSS}. 
\end{rem}

\paragraph{Meaning and reformulation of Theorem  \ref {thm:S-zeta-stokes}.} 

For the convenience of the reader let us  discuss the meaning of
Theorem  \ref {thm:S-zeta-stokes}
 in more detail. 
 Let $V=\FT_\gen(\Fc)_\zeta$,
identified with $\FT_\gen(\Fc)_{\zeta^\pm}$ by monodromy (along the shortest path). 
 Choose $R\gg 0$ and let $\Hen$ be the half-plane $\{ \Re(\zeta \ol w) \geq -R\}$
as in the proof of Proposition \ref {prop:FTgen}, so that
\[
V = H^0_\Hen(\CC, \Fc). 
\]
For any closed subset $Z\subset \Hen$ let $\tau_Z: H^0_Z(\CC,\Fc) \to V$
be the morphism induced by the inclusion of supports.

\begin{figure}[h]
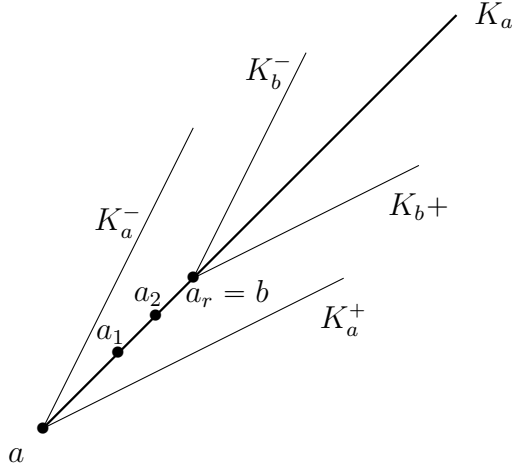

\centering
\btp[scale=0.5]
\node at (0,0){$\bullet$};
\draw (0,0) -- (8,4);

\draw (0,0) -- (4,8);
 
\node at (4,4){$\bullet$};
\node at (2,2){$\bullet$}; 
\node at (1.8, 2.5){$a_1$}; 
\node at (3,3){$\bullet$}; 
\node at (2.8, 3.5) {$a_2$}; 
  
\draw (4,4)-- (10,7);

\draw (4,4)-- (7, 10);

\draw [line width=.3mm] (0,0) -- (4,4);
\draw [line width=.3mm] (4,4) -- (11,11); 
\node at (12,11){$K_a$}; 

\node at (-0.7,-0.7) {$a$};
\node at (4.9, 3.6){$a_r=b$};

\node at (8,2.9){$K_a^+$};
\node at (10, 5.9){$K_b+$};
\node at (2, 5.5){$K_a^-$};
\node at (6, 9.5){$K_b^-$};

\etp
\caption {Crossing a Stokes direction $\zeta=\zeta_{ab}$.}\label{fig:cross-stokes}
\end{figure}

     \vskip .2cm

\vskip .2cm

 As before, for any $a\in\CC$ and $\theta\in S^1$  denote 
$K_a(\theta)=a+\RR_+\theta$ the ray in the direction $\theta$ issuing from $a$.
We denote
\[
K_a^\pm = K_a(\zeta^\pm), \quad K_a = K_a(\zeta), \quad
K^\pm = \bigsqcup_{a\in A} K_a^\pm, \quad K= \bigcup_{a\in A} K_a,
\]
see Fig. \ref {fig:cross-stokes}. 
As $\zeta^\pm$ is non-Stokes,
the rays $K_a^\pm$, $a\in A$, are all distinct, but not so for the $K_a$. 
Let us write $\Phi_a$ for the stalk $\bPhi_a(\Fc)_\theta$ for any $\theta\in [\zeta^+, \zeta^-]$,
these stalks being  identified by  monodromy and put $\Phi=\bigoplus_{a\in A}  \Phi_a$. 
Then we have identifications
\[
\Phi \,  \= \,  \bigoplus_a H^0_{K_a^\pm}(\CC, \Fc) \, =  \, H^0_{K^\pm}(\CC, \Fc). 
\]
Note that the morphisms
\[
\tau_{K^\pm}: \Phi = H^0_{K^\pm}(\CC, \Fc) \to V \text{ as well as } 
\tau_K: H^0_K(\CC, \Fc) \to V
\]
are isomorphisms, since the complements
$\CC\- K^\pm$ and $\CC \- K$ are homotopy
equivalent to $\CC \- \Hen$ and $\Fc$ is locally constant on these complements. 
 It follows from the identifications
constructed in  Proposition \ref {prop:FTgen},
that $S_\zeta= \tau_{K^-}^{-1} \tau_{K^+}$, i.e., it is equal to
 the composition
\[
\Phi = H^0_{K^+} (\CC, \Fc) \buildrel \tau_{K^+} \over  \lra
H^0_{\Hen} (\CC, \Fc)  \buildrel \tau_{K^-}^{-1}  \over \lra 
 H^0_{K^-} (\CC, \Fc)  = \Phi
\]
 Denote for short the inclusion of support maps for individual rays by
 \[
 \tau_{a}^\pm = \tau_{K^\pm(a)}: H^0_{K^\pm_a}(\CC, \Fc)\to V, \quad
 \tau_a = \tau_{K(a)}: H^0_{K(a)}(\CC, \Fc) \to V,
 \]
 so that 
 \[
 \tau_{K^\pm} = \sum_{a\in A} \tau_a^\pm: \,\, \Phi = \bigoplus_{a\in A} \Phi_a 
 \lra  V. 
 \]
 Fix now $a\in A$ and let $a_0=a, a_1, \cdots, a_r =b$ be all elements of $A$
 on $K_a$, see Fig.  \ref{fig:cross-stokes}. Recalling Definition 
 \ref {defi:stokes} of 
  $\St_\zeta$, we can reformulate 
 Theorem  \ref{thm:S-zeta-stokes} as follows.
 
 \begin{refo}\label{refo:Stokes} 
 For any choice of $a$ as above  and any    $\phi\in\Phi_a$ we have
 \[
 \tau_a^-(\phi) \, = \, \tau_a^+(\phi) + \sum_{i=1}^r \tau_{a_i}^+(m_{a, a_i}^-(\phi)). 
 \qed
 \]
 \end{refo} 
 
 \begin{rems}
 (a) Note the similarity of the above formula with Proposition \ref {prop:PL-conseq}(b). 
 Here, instead of the final transport $m_{a_i, b}^+$ to $\Phi_b$,  we have the map $\tau_{a_i}^+$
 to $V$ which corresponds, informally, to putting $b$ at $\oo$. 
 
 \vskip .2cm
 
 (b)  Reformulation \ref{refo:Stokes} matches rather  directly the identities among 
 exponential integrals of multivalued functions along various paths, traditionally used
 in resurgence theory. 
 \end{rems}


 \section{Resurgence theory: convolution algebras in $\ol\Perv(\CC)$} \label{sec:resurg}
 
 \subsection{The general program}\label{subsec:gen-prog}
 
 We now outline  an approach to resurgence  as a program of
 extending and applying   the above elementary  theory to a more general concept of
 perverse sheaves. 
 
  \paragraph { Resurgent perverse sheaves: algebras in the convolution category.}
 We propose to consider perverse sheaves on $\CC_w$ (the Borel plane)
  carrying some algebraic structures
 with respect to the convolution operation $*$. For example, associative (commutative or not)
 algebras, i.e., perverse sheaves $\Ac$ with an operation (i.e., morphism) 
 $\Ac * \Ac \to \Ac$
 satisfying  the associativity and possibly commutativity condition. Or, given such algebra $\Ac$,
 we can consider $\Ac$-modules, i.e., perverse sheaves $\Mc$ with an operation
 $\Ac * \Mc \to \Mc$. Other algebraic structures  can be considered (e.g., Lie algebras). 
 
 \vskip .2cm
 
 An algebraic structure  with respect to convolution defined on a perverse sheaf $\Ac$ would give
 a formal convolution operation on its sections over various domains or on the spaces
 of vanishing cycles (whose intuitive meaning is to describe singularity data of  sections). 
 So various formulas of resurgent analysis  involving convolutions, alien derivatives
 and such could be written intrinsically inside the data associated to  $\Ac$.
 Therefore we propose to call perverse sheaves equipped with such algebra structures
 {\em resurgent perverse sheaves}.\footnote{To be more precise, we suggest to use this term for a generalization of the notion of perverse sheaf discussed in the next subsection.} In various concrete examples  sections of resurgent perverse
 sheaves will be represented by actual resurgent functions in the classical sense. 
 
  \vskip .2cm
  
  As Fourier transform takes convolution to fiberwise multiplication, applying it to
  a resurgent perverse sheaf $\Ac$ would give a local system on $\CC^*_z$ the punctured
  $z$-plane with an  algebra structure  (of the corresponding type)
  in the fibers and with a Stokes structure such that
  the Stokes matrices are isomorphisms of algebras. For example, if $\Ac$ is a commutative
  algebra, the Stokes matrices, being isomorphisms of commutative algebras,
  can be thought of as coordinate changes. This would fit, e.g., into the interpretation of
  cluster transformations as Stokes data for appropriate differential (or rather integral)
  equations, see \cite[\S7]{GMN-4d3d}  and \S \ref{subsec:cluster} below.

 \paragraph{Generalized perverse sheaves and their convolution.}
 In order to realize the above program, we need to generalize the concept of perverse
 sheaves. 
 
 \vskip .2cm
 
 First of all, we need perverse sheaves on $\CC$ whose set of singularities $A$ is
 an arbitrary countable (for example, everywhere dense) subset in $\CC$. A typical
 example of such $A$ is a (free) abelian subgroup in $\CC$ of finite rank $r$; 
 it  cannot be discrete, if $r\geq 3$. Examples like this are inevitable since in the classical case
 (Proposition \ref{prop:F*G})
 the singularities of $\Fc * \Gc$ are typically all the sums of a singularity of $\Fc$ and
 a singularity of $\Gc$. So to have an interesting map $\Ac *\Ac \to \Ac$, the set of
 singularities of $\Ac$ needs to be closed under addition. 
 
  \vskip .2cm

 The fundamental object of study should be
 the localized category $\ol\Perv(\CC,A)$. Its objects should have well-defined
 spaces (more precisey, local systems on $S^1$) of vanishing cycles $\Phi_a$, $a\in A$
 and the transport maps $m_{ab}(\gamma): \Phi_a\to\Phi_b$ for some class of paths
 $\gamma$ joining $a$ and $b$. 
 
 \vskip .2cm
 
 Note that  case of a discrete $A\subset \CC$  may seen to be
  covered by the theory of $\Dc$-modules
 and perverse sheaves in the analytic context. However, already the lifting of $\ol\Perv(\CC,A)$
 back into $\Perv(\CC, A)$ in this context is not obvious, since in the classical case $|A|<\oo$
 the generic stalk of the lifted sheaf is the direct sum of all the $\Phi_a$  
 (see \cite{gelfand-MV,  KKP, fressan-jossen})
 which can
 be infinite-dimensional  for $|A|=\oo$ 
 and so falls outside of the theory of analytic $\Dc$-modules. 
 
  \vskip .2cm
 
 For this and other reasons we need perverse sheaves with possibly infinite-dimensional
 stalks or, more generally, perverse sheaves with values in  a more or less 
 arbitrary abelian category $\Cc$. For example, when $\Cc$ is the category of
 pro-finite-dimensional (= locally linearly compact linearly  topological)
 vector spaces, the dual to the category of all vector spaces,
  this approach  would give (perverse) cosheaves of \cite{2402.07343}. 
 Also, one needs to consider various analytic completions (e.g. of the infinite direct
 sum of the vanishing cycles above), intermediate
 between direct sums and direct products and involving convergence conditions.

 \paragraph{ Lefschetz perverse sheaves in infinite dimensions.}\label{par:lef-inf-dim}
 It is a very appealing idea to generalize the construction of the Lefschetz perverse
 sheaves $\Lc^i_S$ from \S \ref{subsec:lefper} to the case when $X$ is some complex
 function space of ``fields"
 and $S$ is the classical action functional corresponding to some physical theory. 
 
 \vskip .2cm

 Indeed, $\Crit(S)$,  the critical locus of $S$,  is the space of solutions of the classical equations of
 motion; if the problem is set up appropriately,  connected components 
 of $\Crit(S)$ are finite-dimensional.  The behavior of $S$
 in the directions ``transverse''
 to $\Crit(S)$  typically has the form 
 \[
 S(x) = f(x_1,\cdots, x_m) + \sum_{i=m+1}^\oo x_i^2,
 \]
 the direct sum of a function of finitely many variables and an infinite sum of independent squares.
 This means that  perverse sheaf $\Phi = \text{``}\Phi_S(\ul \k _X[\dim X])\text{''}$ on $\Crit(S)$ , or at least,
 on some patches of $\Crit(S)$, can be 
 defined\footnote{
 Strictly speaking, on a compoment $C\subset \Crit(X)$ with $S(C)=a\in\CC$ we should
  define $\Phi_{S-a}$, not $\Phi_S$. 
 } ``by hand'' starting from the $\Phi_f$. As adding an extra  independent square transforms
 vanishing cycles in a known way (Kn\"orrer periodicity), we are lead to a natural gerbe
 (of orientation data) whose trivialization defines $\Phi$ completely.   
  This  by now well
 known procedure   
 is axiomatized using the framework of $(-1)$-shifted symplectic structures \cite{brav}.
 The 
 hypercohomology groups of connected components  $C\subset \Crit(S)$  with coefficients in $\Phi$
 appear in the framework of motivic Donaldson-Thomas  (DT for short) theory. 
 
 \vskip .2cm
 
 However, such a procedure treats different components  independently, 
 as the actual values of $S$ on the components  are ignored or lost. A natural refinement
 of the above data would be 
 perverse sheaves\footnote{or objects of the localized category $\ol\Perv$.}
  $\Lc^i_S$ on $\CC$ whose stalks at $a\in \CC$
 would be $H^i(C_a, \Phi)$, where $C_a$ is the union of components $C\subset \Crit(S)$
 with $S(C)=a$. This additional structure would also provide the transport maps
 $m_{ab}(\gamma)$ between motivic  DT-invariants for a certain class of paths $\gamma$ as long as stability structures are incorporated in our framework (e.g. in the case when  $X$ is the stack of objects of a $3$-dimensional Calabi-Yau category of ``geometric origin\rq\rq{}).

 \paragraph{Unlimited analytic continuation:  the analytic pro-\'etale site.}

 Multivalued analytic functions $f(w)$ on the Borel plane appearing in resurgence theory,
 have the remarkable property of  unlimited analytic continuation 
 which has been made precise using slightly different concepts of
 continuation  ``without cut'' (sans coupure) in  \cite{ecalle}
 or ``without end'' (sans fin)  in \cite{CNP}.  Intuitively, such a formalization needs
 to accomodate two features of the functions in question: 
 
 \begin{itemize}
 \item  [(1)] They possess no natural boundaries, beyond which analytic continuation is not possible
 (such  as the unit circle being the natural boundary for the function $\sum_{n} w^{n^2}$).
 
 \item [(2)] But they can have isolated singularities including
 ramification points that  can accumulate on further and further  sheets of the Riemann
 surface. 
 \end{itemize}
 
 \noindent To explain (2),  any ``branch'' of $f(w)$ is defined over a ``sheet'' obtained
 by removing from $\CC$ a discrete set of cuts emanating from a discrete
 set of ramification points ``visible on this sheet''. But after crossing a cut we arrive on a new sheet
 where $f(w)$ has  a new,  still discrete but possibly  larger set of ramification points etc.
 At the end one can have a seemingly paradoxical outcome that $f$ has a non-discrete,
 e.g., everywhere dense set of singularities (understood as points in $\CC$). 
 
 \vskip .2cm
 
 The features (1) and (2) make one think about the 
  Bhatt-Scholze  theory of the pro-\'etale site \cite{bhatt}. Indeed, (1) suggests some \'etale property
  while going to further and further sheets in (2) resembles some   projective limit procedure. 
  So let us sketch a version
  of this theory in
  the analytic situation. We plan to discuss it in detail in the future.
  
  \vskip .2cm  
 
 Let $X$ be a   complex manifold
 of  dimension $d$.  We can consider on $X$ the 
 {\em analytic Zariski} (or {\em ana-Zariski} for short) topology,
 in which the closed sets are analytic subsets $S\subset X$.
Then the open sets are complements $X\- S$ of such $S$. 
 For example, if $d=1$, then an analytic subset in $X$ is just a discrete
 subset,
 possibly infinite.  Thus an ana-Zariski open subset is a complement
 of a discrete subset.

 \vskip .2cm

 Let now $X$ and $Y$ be complex manifolds of the same dimension $d$. 
 
 \begin{definition}\label{elementary covering}
 A holomorphic map $f: Y\to X$ is called an {\em analytic \'etale}
 ({\em ana-\'etale} for short), if:
 \begin{itemize}
 \item  $f$ is a local biholomorphism, i.e., the differential $df$ 
 is invertible everywhere.
 
 \item  There is an ana-Zariski open sets $Y'\subset Y$, $X'\subset X$
 such that $f(Y') = X'$ and moreover, $f: Y'\to X'$ is an unramified
 covering (with possibly infinite fibers). 
 
 \end{itemize}
 \end{definition}

 \noindent An example is given by the exponential map $\exp: \CC \to \CC$, with $S=\{0\}\subset X=\CC$.

\begin{definition}\label{pro-covering}
A  holomorphic map $f: Y\to X$ is called a {\em pro-ana-\'etale}, if there exist:
\begin{itemize} 
\item A projective system 
\[
Y_0 = X \leftarrow Y_1 \leftarrow Y_2 \leftarrow \cdots
\]
with each arrow being ana-\'etale; 

\item   An injective morphism $\wt f: Y \to \varprojlim Y_i$ whose composition with the projection 
$
\varprojlim Y_i\to Y_0=X$ coincides with $f$. 

\end{itemize} 

\noindent Thus Definition \ref{elementary covering} accounts for the desired feature (1), while
Definition \ref {pro-covering} accounts for (2).

\end{definition}

\begin{exas}
(a) The embedding map $Y=\{|w|<1\} \hra  X = \CC_w$ is a local biholomorphism but not pro-ana-\'etale. Such $Y$ can be seen as the Riemann surface of a function with natural boundary. 

\vskip .2cm

(b) It seems plausible that the classical example of the Riemann surface of the inverse of the hyperelliptic integral discussed in 
 \cite[\S 3] {CNP} can be included into the framework of pro-ana-\'etale theory outlined above. We plan to discuss this as well as more general examples in the future.

\end{exas}

To define constructible and then perverse sheaves with possible non-discrete sets of singularities,
one can follow one of the two paths.

 First, each sheaf $\Fc$  on $X$ in the analytic topology
has the {\em \'etale space} $\wt X_\Fc \to X$ obtained by topologizing the union
of all the stalks $\Fc_x, x\in X$. One can consider sheaves  
whose \'etale spaces  have  maximal Hausdorff parts of their connected components 
satisfying  the property of being  pro-ana-\'etale in the above sense, 
with some constructibility conditions imposed on the sheaves.  

Alternatively, one can consider directly some version of  pro-ana-\'etale Grothendieck 
site  on $X$
and work with sheaves on this site.

 . 
 
 
 \subsection{A finitistic example: COHA of a quiver with potential}\label{subsec:COHA}
 
 Let $Q=(I,E)$ be a finite quiver, with $I$ being the set of vertices and
 $E\to  I\times I$ being the set of oriented edges. 
 For any dimension vector $d = (d_i)_{i\in I}$, $d_i\in \ZZ_+$
 we denote $\Rep_d(Q)$ the stack of complex $d$-dimensional representations $V$
 of $Q$. By definition, such a representation associates to each $i\in I$ a $d_i$-dimensional
 $\CC$-vector space  $V_i$ and to each arrow $e \in E$ between $i$ and $j$  a linear operator 
 $\rho_{e}: V_i\to V_j$
  with no further relations. 
  
  \vskip .2cm

 We denote by $\CC\la Q\ra $ the path algebra of $Q$,
 so representations of $Q$ are the same as left  $\CC\la Q\ra $-modules. 
 Fix a {\em potential}, or a {\em cyclic word} in $\CC\la Q\ra$, i.e., an element 
 $s\in\CC \la Q \ra\bigl/  \bigl[ \CC\la Q\ra, \CC\la Q\ra\bigr]$, 
 the quotient by the commutator subspace (not the ideal generated by commutators). More explicitly, $s$
 is represented as a linear combination of closed edge paths in $Q$. 
 A choice of $s$ gives for any $d$  a regular function
 \[
 S_d: \Rep_d(Q) \lra \CC, \quad V\mapsto \tr_V(s)
 \]
 given by taking the trace. Such functions  are
  additive in the following sense:
  in the {\em induction diagram} of stacks
   \[
  \xymatrix{
  & \{0\to E'\to E \to E''\to 0\} \ar[dl] \ar[dr] & 
  \\
  \Rep_{d'}\times\Rep_{d''} && \Rep_d
  }\begin{matrix} \dim (E')=d',  \dim(E'')=d'',  \\ \dim(E)=d=d'+d''
  \end{matrix}
  \]
  the pullback of $S_{d'} + S_{d''}$ from the left is equal to the pullback
  of $S_d$ from the right. 
  
  \vskip .2cm
 
 The  Cohomological Hall Algebra (COHA) associated to $(Q,s)$
 is \cite{1006.2706}
 \[
A =  \bigoplus_{d\in\ZZ_+^I} A_d, \quad A_d= H^\bullet\bigl(\Rep_d(Q), \Phi_{S_d}(\ul\k[\dim \, \Rep_d(Q)])\bigr). 
 \]
 Here $\dim \,\Rep_d(Q)$ is the dimension in the sense of stacks. 
 The multiplication $A_{d'} \otimes A_{d''}\to A_{d'+d''}$ is given by
 the pullback and pushforward in the induction diagram above.

 Now consider the (bi)graded perverse sheaf
 \[
 \Ac  = \bigoplus_d \Ac_d, \quad \Ac_d = \Lc^\bullet _{S_d}
 \]
 where $\Lc_{S_d}^\bullet$ is the graded  Lefschetz perverse sheaf (see 
 \S \ref{subsec:lefper} \ref{par:lefper}) on $\CC$ associated to $S_d$. 
Each $\Ac_d$ is a graded perverse sheaf with finitely many singularities only. 
We will consider $\Ac_d$ as an object of the localized category $\ol\Perv(\CC)$.
It seems plausible that 
the induction diagram defines morphisms
of graded perverse sheaves  $\Ac_{d'} * \Ac_{d''} \to \Ac_{d''}$ in $\ol\Perv(\CC)$, 
i.e., makes $\Ac$ into an associative 
   convolution algebra refining $A$, cf. \cite[\S4]{1006.2706}.

 \begin{remark}\label{relation to WCS}
 The above definition of COHA depends only on the quiver with potential.
One can make an additional choice consisting  of the central charge $Z:\ZZ^I\to \CC$. This choice gives rise to the wall-crossing structure in the sense of 
\cite{1303.3253, 2005.10651}. If the wall-crossing structure is analytic in the sense of 
\cite{2005.10651}  then it was conjectured in the loc.cit. that germs of sections of the associated non-linear fiber bundle over $\CC$ are resurgent. 
Assuming the conjecture we obtain a perverse sheaf on the Borel plane with  singularities  belonging to the image $Z(\ZZ^I)$. This perverse sheaf seems quite  different from the $\Lc^\bullet_{S_d}$. 
 \end{remark}
 

 \subsection{ Cluster perverse sheaves and wall-crossing structures}\label{subsec:cluster}
 
 \paragraph{Two types of examples: Lefschetz type and cluster type.} 
 Two main classes of perverse sheaves on $\CC$ of infinite rank with possibly infinite set of singularities appear naturally in resurgence theory:
 
 \begin{itemize}
 
 \item  Those which comes from holomorphic functions on infinite-dimensional manifolds (natural generalization of Lefschetz sheaves) 
 
 \item  Those which come from wall-crossing structures  \cite{2005.10651}. 
 In this section we consider this class. 
 They  can be called {\it cluster perverse sheaves\rq\rq{}} or {\it wall-crossing perverse sheaves}. We will use the former term, since cluster transformations play the key role in the story. 
 
 \end{itemize}
 
 These classes  have nonempty intersection. 
 
 \paragraph{Reminder on stability data on graded Lie algebras.} 
 
 Let us recall the relevant structures following  \cite{0811.2435, 1303.3253, 2005.10651}.

 Let $\Gamma$ be a free abelian group of finite rank $n$  endowed
with a skew-symmetric integer-valued  bilinear form
 $\langle - , - \rangle : \Gamma \times \Gamma\to \ZZ$. Consider  the  vector space
$\frak g = \frakg_{\Gamma}=\oplus_{\gamma\in \Gamma}\QQ\cdot e_{\gamma}$.
  This space is made into a Poisson algebra with the commutative (associative) product and
 the Poisson bracket given by
 \be\label{eq:cluster-g}
 \begin{gathered}
 e_{\gamma_1}e_{\gamma_2}=(-1)^{\langle \gamma_1, \gamma_2\rangle}e_{\gamma_1+\gamma_2}\,
 \\
 \{ e_{\gamma_1},e_{\gamma_2}\}=(-1)^{\langle \gamma_1, \gamma_2\rangle}\langle \gamma_1, \gamma_2\rangle e_{\gamma_1+\gamma_2}\,\,.
 \end{gathered}
 \ee
  Let  ${\mathbb T} = \TT_\Gamma:= \Spec( \frakg)$ be the 
  algebraic Poisson manifold
  obtained as  the spectrum of the commutative algebra $\frakg$.
   It  is a torsor over the algebraic torus 
${\Hom}(\Gamma, {\GG}_m)$ and the Poisson structure on ${\mathbb T}$  is invariant with respect to the  torus action.

 \vskip .2cm
 
 Let $\Gamma_\RR=\Gamma\otimes\RR$. For any strictly convex cone
 $C\subset \Gamma_\RR$    we denote 
 \[
 \wh \frakg_C \,= \, \prod_{\gamma\in \Gamma\cap C-\{0\}}
\QQ\cdot e_{\gamma}
  \]
 
 
 the completion of $\frakg$ associated to $C$. It inherits the Poisson algebra
 structure. 

\vskip .2cm

Recall \cite{0811.2435} that   {\em stability data} on $\frakg$
consist of a pair $(Z, a)$, where
  $Z: \Gamma\to\CC$ is a homomorphism
 of abelian groups (``central charge'')  and 
   $a$ is a collection of elements
 $a_\gamma\in \QQ\cdot e_\gamma \subset \frakg$
 given for each  $\gamma\in \Gamma \- \{0\}$. These data satisfy the so-called
  support condition which means that there exists a non-zero quadratic form on $\Gamma\otimes \RR$ which is non-negative on those $\gamma\in \Gamma$ for which $a_\gamma\ne 0$ and which is negative on $\Ker(Z\otimes \RR)$, see \cite[\S2]{0811.2435} for details. 
  
  \vskip .2cm
  
  It is convenient to transform  the  collection $a=(a_\gamma)$ into a collection
 of numbers
 $\Omega = (\Omega(\gamma) \in \QQ)_{\gamma\in\Gamma \- \{0\}}$
 defined uniquely by the identities 
 \[
 a(\gamma)=-\sum_{n\geqslant 1, {1\over{n}}\gamma\in \Gamma\setminus \{0\}}{\Omega(\gamma/n)\over{n^2}}e_{\gamma}\,\,.
 \]
  Then  for any  cone $C$ as above containing $\gamma$ 
  we have a formal identity in $\wh \frakg_C$:
\[
\exp\left(\sum_{n\geqslant 1}a(n\gamma)\right)=
\exp\left(-\sum_{n\geqslant 1} \Omega(n\gamma)\sum_{k\geqslant 1}{e_{kn\gamma}\over{k^2}}\right):=
\exp\left(-\sum_{n\geqslant 1}\Omega(n\gamma){\Li}_2(e_{n\gamma})\right),
\]
where ${\Li}_2(t)=\sum_{k\geqslant 1}t^k/ k^2$ is the dilogarithm series.

\vskip .2cm

As for any Poisson manifold, the
 Lie algebra $\frakg =( \Oc(\TT), \{ -,-\})$ acts on ${\mathbb T}$ by Hamiltonian vector fields (derivations of the coordinate ring $(\frakg = \Oc(\TT), \bullet)$) which have the form
 $\{ f, -\}$ for  $f\in \frakg$. 
  For any $\gamma\in \Gamma\- \{0\}$ 
  denote by $S_{\gamma}$ the formal Poisson automorphism defined by
\[
S_{\gamma}=\exp(\{-{\Li}_2(e_{\gamma}), - \}), \quad  S_\gamma(e_\mu)=(1-e_\gamma)^{\langle \gamma,\mu \rangle} e_\mu
\]
Here the second equality exhibits  $S_\gamma$ as a {\em birational automorphism}
of $\TT$, i.e., an automorphism of the field of fractions of $\frakg$. 
 The first equality shows its formal series expansion in $\wh \frakg_C$ for any
 $C \ni\gamma$ as above. 
It is explained in  \cite{0811.2435}   how   the $\Omega(\gamma)$ are related to enumerative Donaldson-Thomas invariants of $3$-dimensional Calabi-Yau categories.

 \vskip .2cm
 
Further, if the above data satisfy a certain {\it analyticity assumption},  then they give rise to an analytic fiber bundle $E$  over $\CC$ with fibers  isomorphic to
 ${\mathbb T}$.  The gluing functions of $E$ come from transformations
 \be\label{eq:cluster-stokes}
 S_l=\prod_{Z(\gamma)\in l}S_\gamma^{\Omega(\gamma)} 
 \quad \text{ (product in the order given by  $l$)} 
 \ee
 associated to various rays $l\subset \CC$. 
  In this case, the {\it resurgence conjecture} of \cite{2005.10651} 
 says that with any 
  analytic section of $E$  and each $\gamma\in \Gamma$ one can naturally
   associate a resurgent series in the standard coordinate $z$ on $\CC$.
   
   The notions of wall-crossing structure and analytic wall-crossing structure (see \cite{2005.10651}) generalizes the notion of stability data on a graded Lie algebra roughly by considering sheaves of stability data and analytic stability data.

\paragraph{ A perverse sheaf interpretation.} 
From the point of view of the present paper, it is natural to think
of the target of the central charge map $Z:  \Gamma\to\CC$ as the Borel
plane. Let us assume for simplicity that $Z$ is a set-theoreical embedding.
 The transformations $S_l$ are suggestive of Stokes multiplies
and are in fact interpreted as such (in the more general context
of integral equations) in \cite{GMN-4d3d}, cf. also \cite{bridgeland}. 

So it is natural to look 
for a perverse sheaf\footnote{with $\k=\QQ$.} $\Fc$  on $\CC$ in some generalized sense as above
(in fact, an object of an apprioriate
localized category $\ol\Perv$)  with the properties: 

\begin{itemize}
\item [(0)] $\Fc$ is a Poisson algebra with respect to convolution. 

\item[(1)] The set $A$ of singularities of $\Fc$ is a subgroup of $Z(\Gamma)$. 
So it can be discrete (though still infinite), if $\rk (\Gamma) \leq 2$.

\item[(2)] For $a=Z(\gamma)$ the space $\Phi_a(\Fc)$ is $1$-dimensional,
identified with $\QQ\cdot e_\gamma$. 

 \item[(3)] For $\zeta\in S^1$ the  Stokes operator $\St_\zeta$ for $\Fc$
is given by $S_{\RR_+\zeta}$ from \eqref {eq:cluster-stokes}. 

\end{itemize}

Then $\frakg$  would be realized as  $\bigoplus_{a\in A} \Phi_a(\Fc)$ with the operations coming
from those on $\Fc$. The Fourier transform $\wh\Fc$ would  then
be a generalized perverse sheaf of infinite rank (the latter needs a rigorous definition),
 and the stalks of the local system 
 $\FT(\Fc)_\gen$ on $\CC^*$ would involve some analytic completions of 
 $\frakg = 
 \bigoplus_{\gamma\in\Gamma}  \QQ\cdot e_\gamma$.
 For a ray $l\subset \CC$ the transformation 
   $S_l$ need not preserve $\frakg$  which is an algebraic
   direct sum: it maps $\frakg$  into the direct product. 
    But as such  it has a well defined matrix element between any two summands
    which we can write as $(S_l)_{ab}: \Phi_a(\Fc)\to\Phi_b(\Fc)$. 
    Note that $(S_l)_{ab} = 0$ unless $(b-a)\in l$. Therefore it is  indeed meaningful
    to interpret each $S_l$ as the Stokes operator of some would-be
    generalized perverse sheaf $\Fc$. If we do so, 
    the finite case analysis
   of \S \ref{subsec:alien-der} provides an answer for what should be
   the (say, left-)avoiding transport  $m^{-, \Fc}_{ab}$ or the alien transport
 $m^{\Delta, \Fc}_{ab}$ for any $a,b\in A= Z(\Gamma)$. 
 Proposition \ref{prop:xi-delta} suggests that his should be
 sufficient information to recover $\Fc$ as an object of the localized category. 
 
 \vskip .2cm

 This picture  is expecially compelling when $\rk(\Gamma)=2$ and $A$ is a
 discrete lattice in $\CC$. Then one can look for $\Fc$ as a perverse sheaf in the
 classical sense but with possibly infinite-dimensional stalks.

 
 \subsection{Perverse sheaves in Chern-Simons theory} \label{subsec: CS}
 An example of   a holomorphic function on an infinite-dimensional complex manifold (or stack)
 is provided by the 
  complexified Chern-Simons functional (CS functional for short)
   associated to a compact oriented $3$-manifold $M$
  and a complex semisimple Lie group $G$. Let e.g. $G=SL_n(\CC)$. 
  We take  $\Xc$  to be the moduli stack
  of $C^\oo$-connections on the trivial $SL_n$-bundle on $M$, so $\Xc$ is the quotient
  of  the vector space $\Omega^1(M)\otimes \sen\len_n$ by the gauge group 
  $C^\oo(M, SL_n(\CC))$, and the CS functional is
  \be\label{eq:CS}
  \CS: \Xc \lra \CC/4\pi^2 \ZZ, \quad \CS(A) = \int_{M^3} \tr \left({1\over{2}}A\wedge dA+{1\over{3}}A\wedge A\wedge A\right)\, \on{mod}\,  4\pi^2\ZZ. 
  \ee
  In fact, quotienting by the subgroup of gauge transformations equal to $1$ over a fixed point
  $m_0\in M$, gives  an infinite-dimensional {\em manifold} $\Xc_\fr$ (framed connections) with
$\Xc = \Xc_\fr / SL_n(\CC)$  and we can consider $\CS$ on this manifold.

\vskip .2cm

In order  to  get rid of the  multivaluedness we pass to the maximal abelian covering
$\wt \Xc_\fr \buildrel \ZZ \over \to  \Xc_\fr$, so we get a well defined holomorphic function
\be\label{eq:tildeCS}
\wt\CS: \wt\Xc_\fr \lra \CC.
\ee
 One can   then attempt to define and study the corresponding Lefschetz perverse sheaf $\Lc_\CS$
  along the lines discussed in \S \ref{sec:resurg}\ref{par:lef-inf-dim}. 
  Such a study was  in fact initiated in  \cite  {2402.07343} although 
  very little is known about $\Lc_{\wt\CS}$. Some known facts and some expectations can be summarized as follows.
  
  \begin{itemize}
  \item[(1)] The  functional $\CS$ in \eqref{eq:CS} has finitely many critical values,
  which are known  to be the regulators of some elements in algebraic K-group $K_3(\CC)$.
  Accordingly, the critical values of $\wt\CS$ fall into finitely many arithmetic 
  progressions with step $4\pi^2\ZZ$. 
  
  \item[(2)] The generic stalks of the expected perverse sheaf $\Lc_{\wt\CS}$
  (i.e., intuitively, the middle-dimensional cohomology of the generic fiber of $\wt\CS$)
   are
  infinite-dimensional. In fact, it is easier to  first define the Verdier dual cosheaf, as done in
  \cite  {2402.07343}. But the spaces of vanishing cycles are finite-dimensional. 
  
  \item[(3)] $\Lc_{\wt\CS}$ can be defined as the perverse extension of the
  local system formed by the middle cohomology of the fibers
   on the  complement to the set of critical values,
  see \cite[\S8.3] {2402.07343}. 
  
  \item[(4)] In addition,  $R\Gamma(\CC, \Lc_{\wt\CS})=0$, so $\Lc_{\wt\CS}$ is, formally,
  an object of the category $\Perv^0(\CC)$ (although with infinitely many singularities
  and  with infinite-dimensional stalks). 
  
  \end{itemize}

 \noindent  Some further conjectures can be found in \cite [\S 8.3]{2402.07343}. 
   It was also explained in {\em  loc.cit.}  how the  wall-crossing structure and resurgence properties of the perturbative expansions of the Chern-Simons functional integral are related to the perverse sheaf
 $\Lc_{\wt\CS}$.

  \vskip .2cm

    \vfill\eject

\vfill\eject

	\vskip.5cm

M. K.:  Kavli IPMU, 5-1-5 Kashiwanoha, Kashiwa, Chiba, 277-8583 Japan. 
\\
  { mikhail.kapranov@protonmail.com}
  
  \vskip .2cm

Y.S.: { Department of Mathematics, Kansas State University, Manhattan KS 66508 USA.
{ soibel@math.ksu.edu} 

} 

\end{document}